\documentclass[12pt]{article}
\usepackage{amsmath}
\usepackage{amssymb}
\newsymbol \blackbox 1004
\newcommand{\eh}{\hfill}\newlength{\sperr}

\newenvironment{proof}{{\settowidth{\sperr}{\bf\rm
Proof}%
\par\addvspace{0.3cm}\noindent\parbox[t]{1.3\sperr}
{\bf\rm P\eh r\eh o\eh o\eh f\eh }%
}}{\nopagebreak\mbox{} $\blackbox$\par\addvspace{0.3cm}}

\def\a{\alpha}
\def\b{\beta}

\def\s{\sigma}

\def\t{\theta}

\def\vp{\varphi}
\def\ve{\varepsilon}
\def\wh{\widehat}
\def\wt{\widetilde}
\def\ov{\overline}

\def\BC{{\mathbb C}}
\def\BR{{\mathbb R}}

\newtheorem{Pa}{Paper}[section]
\newtheorem{Tm}[Pa]{{\bf Theorem}}

\newtheorem{Cy}[Pa]{{\bf Corollary}}
\newtheorem{Rk}[Pa]{{\bf Remark}}

\newtheorem{Ee}[Pa]{{\bf Example}}
\newtheorem{Dn}[Pa]{{\bf Definition}}
\newtheorem{Pn}[Pa]{{\bf Proposition}}
\newcommand{\CC}
{{\mathchoice {\setbox0=\hbox{$\displaystyle\rm C$}\hbox{\hbox
to0pt{\kern0.4\wd0\vrule height0.9\ht0\hss}\box0}}
{\setbox0=\hbox{$\textstyle\rm C$}\hbox{\hbox
to0pt{\kern0.4\wd0\vrule height0.9\ht0\hss}\box0}}
{\setbox0=\hbox{$\scriptstyle\rm C$}\hbox{\hbox
to0pt{\kern0.4\wd0\vrule height0.9\ht0\hss}\box0}}
{\setbox0=\hbox{$\scriptscriptstyle\rm C$}\hbox{\hbox
to0pt{\kern0.4\wd0\vrule height0.9\ht0\hss}\box0}}}}

\title{Weyl matrix functions and inverse problems for discrete Dirac type
self-adjoint system: explicit and general solutions}

\author{B. Fritzsche, B. Kirstein, I.Ya. Roitberg,  A.L. Sakhnovich}

\date{}
\begin{document}
\maketitle

\begin{abstract}
Discrete Dirac type self-adjoint system is equivalent to the block
Szeg\"o recurrence. Representation of the fundamental solution is
obtained, inverse problems on the interval and semiaxis are
solved. A Borg-Marchenko type result is obtained too. Connections
with the block Toeplitz matrices are treated.
\end{abstract}

{\bf MSC2000: 39A12, 37K35, 47B35}

{\it Keywords: Discrete Dirac system, Szeg\"o recurrence, Weyl
function, inverse problem, $j$-theory, block Toeplitz matrix}

\section{Introduction} \label{intro}
 \setcounter{equation}{0}
Continuous self-adjoint Dirac type system
\begin{equation} \label{0.3}
  \frac{d Y}{dx}(x,z)= i\big( z j + j V(x)\big) Y(x, z ),
  \quad
j = \left[
\begin{array}{cc}
I_{p} & 0 \\ 0 & -I_{p}
\end{array}
\right], \, \, V=\left[
\begin{array}{lr}
0 & v \\ v^* & 0
\end{array}
\right]
\end{equation}
is a classical object of analysis with various applications (in
mathematical physics and nonlinear integrable equations, in
particular). Here $I_p$ is the $p \times p $ identity matrix and
$v$ is a $p \times p$ matrix function. In this paper, we treat a
discrete self-adjoint Dirac type system:
\begin{equation} \label{0.1}
W_{k+1}( \lambda)-W_k(\lambda)=- \frac{i}{\lambda}j C_k
W_k(\lambda) \quad (k \geq 0),
\end{equation}
where $C_k$ are $m \times m$ matrices, $m=2p$,
\begin{equation} \label{0.2}
C_k=C_k^*, \quad C_k j C_k=j,
\end{equation}

To see that (\ref{0.1}) is a discrete analog of system
(\ref{0.3}), notice that  (\ref{0.3}) is equivalent to a subclass
of canonical systems $W_x=i z j H(x)W$ (see \cite{SaA7, SaL3} and
references therein). One can follow also the arguments from
\cite{KaS}, where the skew self-adjoint discrete  Dirac type
system  have been studied and explicit solutions of the isotropic
Heisenberg magnet model have been obtained. As suggested in
\cite{KaS} introduce matrix functions $U$ and $W$ by the relations
\begin{equation} \label{0.4}
W(x,z)= U(x)Y(x,z), \quad \frac{d U}{dx}(x)=-i U(x)j V(x), \quad
U(0)=I_{m}.
\end{equation}
Since $V$ is self-adjoint, we get from (\ref{0.4}) that $U$ is
$j$-unitary, i.e., $UjU^* \equiv j$. Now (\ref{0.3}) and the first
relation in (\ref{0.4}) yield
\begin{equation} \label{0.5}
  \frac{d W}{dx}=  \frac{d U}{dx}U^{-1}W+iU \big( z j + j V\big) U^{-1}W=izjHW,
\end{equation}
where $H=jUU^*j=H^*$, $HjH \equiv j$. Compare system (\ref{0.1}),
where matrices $C_k$ satisfy (\ref{0.2}), and system (\ref{0.5})
to see that (\ref{0.1}) is an immediate discrete analog of
(\ref{0.3}).

When $p=1$ and $C_k>0$, system (\ref{0.1}) is equivalent to the
well-known self-adjoint Szeg\"o recurrence, which plays an
important role in the orthogonal polynomials theory and  is also
an auxiliary system for the Ablowitz-Ladik hierarchy (see, for
instance, \cite{GGH, GHMT, Si2} and various references therein).
The equivalence of system (\ref{0.1}), where $C_k>0$ and $C_k j
C_k=j$, to the block (matrix-valued) Szeg\"o recurrence  is given
in Proposition \ref{equiv}.

We consider representation of the fundamental solution of system
(\ref{0.1}) and solve direct and inverse problems directly in
terms of the Weyl functions. Both explicit and general solutions
are obtained. First, we  obtain explicit solutions of the direct
and inverse problems for system (\ref{0.1})
 for the case of the so called
{\it pseudo-exponential potentials} $C_k$ (the case of the
rational Weyl functions). Our case includes as a subcase the
rapidly decaying {\it strictly pseudo-exponential potentials}.
Recall that discrete and continuous systems with the potentials,
which belong to the subclass of the  strictly pseudo-exponential
potentials, have been actively studied in \cite{AG0}-\cite{AG2},
\cite{AGKS}-\cite{AD}. In particular, direct and inverse problems
for Szeg\"o recurrence on the semiaxis with the scalar ($p=1$)
strictly pseudo-exponential potentials have been treated in
\cite{AG1, AG2}. Direct and inverse problems for the
pseudo-exponential potentials (continuous case) have been studied
in a series of Gohberg-Kaashoek-Sakhnovich papers \cite{GKS1,
GKS6} (see references therein and see also \cite{FKS} for the case
of the generalized pseudo-exponential potentials). The case of the
discrete skew-self-adjoint Dirac system have been studied in
\cite{KaS}. Notice that similar to \cite{GKS1, SaA5} (see also
\cite{AG1, AG2, AGKS, GKS1, KaS}) we start our explicit
constructions with the explicit formula for the fundamental
solution.

For a more general (non-rational) situation of the Weyl functions
$\displaystyle \phi(z)=\sum_{k=0}^\infty \phi_k z^k$ such that
$\sum_{k=0}^\infty  \|\phi_k \|<\infty$ the direct  problem for a
block (matrix-valued) Szeg\"o recurrence on the semiaxis
(including non-self-adjoint case and under some  additional
conditions) is treated in a recent important Alpay-Gohberg paper
\cite{AG3}. In Sections \ref{gencase} and  \ref{Toepl} we solve
direct and inverse problems for the general type potentials
$C_k>0$ (and thus for the general type self-adjoint block Szeg\"o
recurrence) on the interval and semiaxis. Borg-Marchenko type
uniqueness result for system (\ref{0.1}) is obtained too.
Connections with the well known Toeplitz matrices appear. For the
interesting discussions on the connections between Toepliz
matrices, Szeg\"o recurrencies and orthogonal polynomials see also
\cite{AGtep, Dy} and references therein. Interesting spectral
theoretical results on the discrete canonical systems, where
$C_kjC_k=0$, one can find in \cite{RoSa, SaL2}. A complete Weyl
theory for the Jacobi matrices and various useful references are
contained in \cite{GT}.

\section{Preliminaries} \label{Prel}
 \setcounter{equation}{0}
An important discrete analog of Dirac type system takes the form
\begin{equation} \label{0.6}
X_{k+1}(z)=\t_k R_k \left[
\begin{array}{cc}
z I_{p} & 0 \\ 0 & I_{p}
\end{array}
\right] X_k(z), \quad R_k=R_k^*, \quad R_k j R_k^*=j
\end{equation}
$(\t_k \in \BC, \, \t_k \not=0)$.  System  (\ref{0.6}) can be
presented in the form (\ref{0.1}) after transformation
\begin{equation} \label{0.7}
W_k(\lambda)=\frac{\big(i- \lambda^{-1} \big)^k}{\prod_{r=0}^{k-1}
\t_r} U_k \left[
\begin{array}{cc}
z I_{p} & 0 \\ 0 & I_{p}
\end{array}
\right]X_k(z), \quad C_k=(U_k^*)^{-1}R_k^2U_k^{-1},
\end{equation}
where $U_0:=I_m$,
\begin{equation} \label{0.8}
U_k:=(ij R_0)(ijR_1)\times \ldots \times (ij R_{k-1}) \quad (k>0),
\quad z= \frac{1+i \lambda }{1-i \lambda }.
\end{equation}
 A particular scalar case
($p=1$) of system (\ref{0.6}) is a well known Szeg\"o recurrence,
where
\begin{equation} \label{n.1}
R_k=\frac{1}{\sqrt{1-|\rho_k|^2}} \left[
\begin{array}{cc}
1 & - \rho_k \\ - \ov{ \rho_k} & 1
\end{array}
\right]>0, \quad |\rho_k|<1, \quad \t_k=\sqrt{1-|\rho_k|^2}.
\end{equation}
Further we assume that some expression for $\theta_k$ via $R_k$ is
fixed for $p>1$ too. For $p=1$ representations (\ref{n.1}) of
$R_k$ follow from the relations $R_k=R_k^*>0$ and $R_k j R_k=j$.
Coefficients $\rho_k$ are called Schur (or sometimes Verblunsky)
coefficients (see, for instance, \cite{DFK, Si2} and various
references therein). Notice that matrices $C_k$ given by the
second relation in (\ref{0.7}) are positive definite. Vice versa,
if  matrices $C_k$
 are positive definite, Szeg\"o recurrence is uniquely recovered
from  system (\ref{0.1}). The same is true for the block Szeg\"o
recurrences (\ref{0.6}), where $R_k>0$ and $\theta_k=\theta(R_k)$
for some function $\t$.
\begin{Pn} \label{equiv}
There is a one to one correspondence between the subclass of
system (\ref{0.1}), where the matrices $C_k>0$ satisfy
(\ref{0.2}), and block Szeg\"o recurrences (\ref{0.6}), where
$R_k>0$ and $\theta_k=\theta(R_k)$. This correspondence is given
by (\ref{0.7}),  (\ref{0.8}) to map block Szeg\"o recurrences into
Dirac type systems and is given by the equality
\begin{equation} \label{n.2}
U_0=I_m, \quad R_k=\Big( U_k^*C_kU_k \Big)^\frac{1}{2}>0, \quad
U_k=U_{k-1}(i j R_{k-1})
\end{equation}
to map Dirac type systems into block Szeg\"o recurrences.
\end{Pn}
\begin{proof}.
The first part of the proposition is already proved. Moreover,
according to (\ref{0.7}) and (\ref{0.8}) matrices $R_k$ ($k \geq
0$) and $U_k$ ($k>0$) are uniquely defined by the relations
(\ref{n.2}). Clearly we have $R_k=R_k^*>0$,
 and it remains to
prove $R_k j R_k=j$. We shall use for this purpose a unitary
equivalence of $ U_k^*C_kU_k$ to diagonal matrix $D_k>0$:
\begin{equation} \label{n.3}
U_k^*C_kU_k=\wh U_k^*D_k\wh U_k, \quad \wh U_k^* \wh U_k=\wh U_k
\wh U_k^*=I_m.
\end{equation}
When $k>0$ we assume that $R_r j R_r=j$ is already proved for
$r<k$, and so $U_k j U_k^*=j$. Then formula (\ref{n.3}) implies
that $\wh U_k^*D_k\wh U_k j \wh U_k^*D_k\wh U_k=j$, i.e.,
\begin{equation} \label{n.4}
D_k^{-1}=J_k D_k J_k, \quad J_k:=\wh U_k j \wh U_k^*, \quad
J_k=J_k^*=J_k^{-1}.
\end{equation}
By (\ref{n.4}), without loss of generality we can assume that the
diagonal matrix $D_k$ has the following form:
\begin{equation} \label{nn.1}
{\displaystyle{D_k={\mathrm{diag}}\{d_0 I_{l_0}, \, d_1 I_{l_1},
\, d_1^{-1} I_{l_1}, \ldots, \, d_r I_{l_r}, \, d_r^{-1}
I_{l_r}\}}},
\end{equation}
where $d_0=1$ and $d_l \not= d_s$, $d_l \not= d_s^{-1}$ for $l
\not= s$. Formulae (\ref{n.4}) and (\ref{nn.1}) imply, in their
turn, that the matrix $J_k$ has a block diagonal structure
\begin{equation} \label{nn.2}
{\displaystyle{J_k={\mathrm{diag}}\{u_0, \, j_1, \ldots, \, j_r
\}}}, \quad j_s=\left[
\begin{array}{cc}
0 & u_s \\ u_s^* & 0
\end{array}
\right],
\end{equation}
where $u_s$ are $l_s \times l_s$ unitary matrices, and
$u_0=u_0^*$. According to (\ref{nn.1}) and (\ref{nn.2}) we have
\begin{equation} \label{n.5}
J_kD_k^\frac{1}{2}=D_k^{-\frac{1}{2}}J_k.
\end{equation}
By  (\ref{n.3}) we have $R_k=\wh U_k^*D_k^\frac{1}{2}\wh U_k$,
and, taking into account
 (\ref{n.5}), we  get
\[
R_k j R_k=\wh U_k^*D_k^\frac{1}{2} J_kD_k^\frac{1}{2}\wh U_k=\wh
U_k^* J_k \wh U_k= \wh U_k^* \wh U_k j \wh U_k^* \wh U_k=j.
\]
\end{proof}

The spectral theory of the discrete and continuous systems is
strongly related to the construction of the fundamental solutions
(see, for instance, \cite{AG2}-\cite{AD}, \cite{GKS1, GKS6, KaS,
Po}, \cite{SaA5}-\cite{SaL3} and references therein). The
$j$-properties of the fundamental solutions play an important role
\cite{AD,  DFK, FK, FKMo, FKMu, GoNe, Po, SaL2, SaL3}.

For the case of the explicit construction the version of the
Backlund-Darboux transformation (BDT) introduced in \cite{SaA1,
SaA2, SaA5} proves very fruitful. Choose $n>0$, two $n \times n$
parameter matrices $A$ ($\det A \not=0$) and $S_0=S_0^*$, and $n
\times m$ parameter matrix $\Pi_0$ such that
\begin{equation} \label{0.9}
A S_0-S_0 A^*=i \Pi_0 j \Pi_0^*.
\end{equation}
Define sequences $\{\Pi_k\}$ and $\{S_k\}$  ($k >0$) by the
relations
\begin{equation} \label{0.10}
\Pi_{k+1}=\Pi_k+i A^{-1}\Pi_k j,
\end{equation}
\begin{equation} \label{0.11}
S_{k+1}=S_k+A^{-1}S_k (A^*)^{-1}+A^{-1}\Pi_k \Pi_k^*(A^*)^{-1}.
\end{equation}
It follows that the matrix identity
\begin{equation} \label{0.11'}
A S_{k+1}-S_{k+1} A^*=i \Pi_{k+1} j \Pi_{k+1}^* \quad (k \geq 0)
\end{equation}
is true.  Following the lines of the discrete BDT version for the
skew self-adjoint discrete  Dirac type system presented in
\cite{KaS}, we get the  theorem.
\begin{Tm} \label{FundSol} Suppose $\det \, S_r \not=0$ $(0 \leq
r \leq N)$. Then the fundamental solution $W_{k+1}$ of system
(\ref{0.1}), where
\begin{equation} \label{0.13}
C_k:=I_m+\Pi_k^*S_k^{-1}\Pi_k-\Pi_{k+1}^*S_{k+1}^{-1}\Pi_{k+1},
\end{equation}
admits representation
\begin{equation} \label{0.15}
W_{k+1}(\lambda)=w_A(k+1,\lambda)\big(I_m - \frac{i}{\lambda} j
\big)^{k+1}w_A(0,\lambda)^{-1} \quad (0 \leq k < N).
\end{equation}
Here $W_{k+1}$ is normalized by the condition $W_0(\lambda)=I_m$,
and
\begin{equation} \label{0.14}
w_A(k,\lambda):=I_m-i j \Pi_k^*S_k^{-1}(A-\lambda I_n)^{-1}\Pi_k,
\end{equation}
\end{Tm}
The right hand side of  (\ref{0.14}) with fixed $k$ is a so called
transfer matrix function in Lev Sakhnovich form
\cite{SaL1}-\cite{SaL3}.

We say that system (\ref{0.1}), where matrices $C_k$ are given by
(\ref{0.13}), is {\it determined} by the parameter matrices $A$,
$S_0$ and $\Pi_0$.

\begin{proof} of theorem.
Formula (\ref{0.15}) easily follows from the basic for this proof
equality
\begin{equation} \label{0.12}
w_A(k+1, \lambda)\big(I_m - \frac{i}{\lambda} j \big)=\big(I_m -
\frac{i}{\lambda} j C_k \big)w_A(k, \lambda),
\end{equation}
that we shall derive now. Taking into account  (\ref{0.14}) one
can see that (\ref{0.12}) is equivalent to the equality
\[
-\frac{i}{\lambda}j (I_m - C_k)=- \big(I_m - \frac{i}{\lambda} j
C_k \big)i j \Pi_k^*S_k^{-1}(A-\lambda I_n)^{-1}\Pi_k
\]
\begin{equation} \label{0.16}
+i j \Pi_{k+1}^*S_{k+1}^{-1}(A-\lambda I_n)^{-1}\Pi_{k+1}\big(I_m
- \frac{i}{\lambda} j \big),
\end{equation}
i.e., the Taylor coefficients at infinity of the matrix functions
in both sides of (\ref{0.16}) coincide. Hence, by the series
expansion $\displaystyle (A-\lambda
I_n)^{-1}=-\lambda^{-1}\sum_{r=0}^\infty\big(\lambda^{-1}A\big)^r$
and formula (\ref{0.10}), formula (\ref{0.16}) is equivalent to a
family of equalities:
\begin{equation} \label{0.17}
I_m-C_k=-\Pi_k^*S_k^{-1}\Pi_k+\Pi_{k+1}^*S_{k+1}^{-1}\Pi_{k+1}
\end{equation}
and
\begin{equation} \label{0.18}
K_kA^{r-2}\Pi_k=0 \quad (r>0),
\end{equation}
where
\begin{equation} \label{0.19}
K_k:=\Pi_{k+1}^*S_{k+1}^{-1}(A^2+I_n)-\Pi_{k}^*S_{k}^{-1}A^2+i C_k
j  \Pi_{k}^*S_{k}^{-1}A.
\end{equation}
Notice that (\ref{0.17}) is immediate from (\ref{0.13}). If we
prove also $K_k=0$, then (\ref{0.18}) will follow, and so we will
get (\ref{0.16}) or equivalently (\ref{0.12}), which implies
(\ref{0.15}). It remains to show that $K_k=0$. For this purpose we
shall rewrite (\ref{0.19}) using (\ref{0.10}) and (\ref{0.13}):
\[
K_k=\Pi_{k+1}^*S_{k+1}^{-1}(A^2+I_n)-\Pi_{k}^*S_{k}^{-1}A^2+i  j
\Pi_{k}^*S_{k}^{-1}A+i \Pi_k^*S_k^{-1}\Pi_k j \Pi_{k}^*S_{k}^{-1}A
\]
\begin{equation} \label{0.20}
-i \Pi_{k+1}^*S_{k+1}^{-1}\big( \Pi_k+i A^{-1}\Pi_k j    \big)j
\Pi_{k}^*S_{k}^{-1}A.
\end{equation}
According to (\ref{0.11'}) we have $i \Pi_k j
\Pi_{k}^*S_{k}^{-1}=A-S_kA^*S_{k}^{-1}$. Therefore, from
(\ref{0.20}) we derive
\[
K_k=\Pi_{k+1}^*S_{k+1}^{-1}\big(I_n+S_kA^*S_k^{-1}A+A^{-1}\Pi_k\Pi_k^*S_k^{-1}A
\big)- \Pi_{k}^*A^*S_{k}^{-1}A+i  j \Pi_{k}^*S_{k}^{-1}A.
\]
In view of (\ref{0.11}) we simplify our last formula:
\begin{equation} \label{0.21}
K_k=\Pi_{k+1}^*A^*S_k^{-1}A- \Pi_{k}^*A^*S_{k}^{-1}A+i  j
\Pi_{k}^*S_{k}^{-1}A.
\end{equation}
Finally, by (\ref{0.10}) and (\ref{0.21}) we have $K_k=0$.
\end{proof}
\begin{Pn}\label{Pn2.2}
Suppose $\det \, S_r \not=0$ $(0 \leq r \leq N)$. Then the
matrices $C_k$  $(0 \leq k < N)$ given by (\ref{0.13}) satisfy
conditions (\ref{0.2}).
\end{Pn}

\begin{proof} .
The first equality  in (\ref{0.2}) is immediate. To prove the
second equality notice that by the standard in the $S$-node theory
\cite{SaL1}-\cite{SaL3} calculations (see also, for instance,
formula (2.10) in \cite {FKS} it follows from (\ref{0.9}) and
(\ref{0.11})  that
\begin{equation} \label{0.22}
w_A(r,  \lambda)^*j w_A(r,  \lambda)=j+i(\ov \lambda
-\lambda)\Pi_r^*(A^*-\ov \lambda I_n)^{-1}S_r^{-1} (A- \lambda
I_n)^{-1}\Pi_r.
\end{equation}
In particular, we have
\begin{equation} \label{0.23}
w_A(r,  \ov \lambda)^*j w_A(r,  \lambda)=j. \quad r \geq 0.
\end{equation}
It is easily checked also that
\begin{equation} \label{0.24}
\Big( I_m + \frac{i}{\lambda} j\Big)j\Big( I_m -\frac{i}{\lambda}
j\Big)=\Big( 1 + \frac{1}{\lambda^2} \Big)j.
\end{equation}
According to (\ref{0.12}) formulas (\ref{0.23}) and  (\ref{0.24})
yield the equality
\begin{equation} \label{0.25}
\Big( I_m + \frac{i}{\lambda} C_k j\Big)j\Big( I_m
-\frac{i}{\lambda} j C_k \Big)=\Big( 1 + \frac{1}{\lambda^2}
\Big)j.
\end{equation}
Therefore the second  equality  in (\ref{0.2}) holds.
\end{proof}
\section{Auxiliary propositions} \label{Aux}
 \setcounter{equation}{0}
Recall that the invertibility of matrices $S_k$ is essential for
our constructions. On the other hand the important subcase of
Szeg\"o recursion corresponds to system (\ref{0.1}), where
$C_k>0$. A natural condition, when all $S_k>0$ and $C_k>0$ is
given in our next proposition.
\begin{Pn} \label{Pn3.1}
Let the parameter matrix $S_0$ be positive definite, i.e.,
$S_0>0$. Then we have
\begin{equation} \label{1.1}
S_k>0 \quad (k\geq 0), \quad C_k>0 \quad (k\geq 0).
\end{equation}
\end{Pn}
\begin{proof}.
The inequalities for $S_k$ in  (\ref{1.1}) follow from
(\ref{0.11}) by induction. To derive the relations $C_k>0$,
introduce first two block matrices:
\begin{equation} \label{1.2}
G=\left[
\begin{array}{cc}
S_k & \Pi_k \\ \Pi_k^* & c I_{m}
\end{array}
\right], F=\left[
\begin{array}{cc}
A^{-1} & a A^{-1} \Pi_k \\ 0 & -i b j
\end{array}
\right],
\end{equation}
where
\begin{equation} \label{1.2'}
a(2+ac)=1, \quad b(1+ac)=1,
\end{equation}
 and, moreover, $c$ is sufficiently large so that
$G>0$.  We shall discuss the choice of $a$ and $b$ satisfying
(\ref{1.2'}) later on. According to (\ref{0.10}), (\ref{0.11}),
(\ref{1.2}) and (\ref{1.2'}), direct calculations show that
\begin{equation} \label{1.3}
G+FGF^*= \left[
\begin{array}{cc}
S_{k+1} & \Pi_{k+1} \\ \Pi_{k+1}^* & c(1+b^2) I_{m}
\end{array}
\right].
\end{equation}
As $G+FGF^*>G>0$ we have $G^{-1}>(G+FGF^*)^{-1}$, and, therefore,
the inequality holds also for the $m \times m$ right lower blocks
of these matrices:
$\big(G^{-1}\big)_{22}>\big((G+FGF^*)^{-1}\big)_{22}$. Finally, we
obtain
\begin{equation} \label{1.4}
\Big(\big(G^{-1}\big)_{22}\Big)^{-1}<\Big(\big((G+FGF^*)^{-1}\big)_{22}\Big)^{-1}.
\end{equation}
Taking into account (\ref{1.2}), we can rewrite (\ref{1.4}) in the
form
\begin{equation} \label{1.5}
cI_m-\Pi_k^*S_k^{-1}\Pi_k<c(1+b^2)I_m -
\Pi_{k+1}^*S_{k+1}^{-1}\Pi_{k+1}.
\end{equation}
Let us fix $c$ and choose a root $a$ $(0<a<1/2)$, of the equation
\[
a^2+\frac{2}{c}a-\frac{1}{c}=0,
\]
which is always possible. Putting also $b=a(1-a)^{-1}$, we see
that relations (\ref{1.2'}) hold. Moreover, the first relation in
(\ref{1.2'}) means that $a^2c=1-2 a$, and so
\begin{equation} \label{1.6}
cb^2=ca^2(1-a)^2=(1-2a)(1-a)^2<1.
\end{equation}
From (\ref{1.5}) and (\ref{1.6}) it follows that
\begin{equation} \label{1.7}
I_m+\Pi_k^*S_k^{-1}\Pi_k- \Pi_{k+1}^*S_{k+1}^{-1}\Pi_{k+1}>0
\end{equation}
Recall the definition (\ref{0.13}) to see that inequality
(\ref{1.7}) implies $C_k>0$.
\end{proof}
In this section we shall need as well another property of $C_k$.
\begin{Pn} \label{Pn3.2}
Let relations (\ref{0.2}) hold, and assume that $C_k>0$. Then we
have $C_k \pm j \geq 0$.
\end{Pn}
\begin{proof}.
It follows from (\ref{0.2})  that
\begin{equation} \label{1.8}
(C_k + \ve j)j(C_k+ \ve j)=2 \ve \big(C_k+ \frac{1+ \ve^2}{2 \ve}
j \big).
\end{equation}
If $(C_k + \ve j) f=0$, then by (\ref{1.8}) we have also $
\big(C_k+ \frac{1+ \ve^2}{2 \ve} j \big)f=0$, and so
$(1-\ve^2)jf=0$. Therefore, we have $\det (C_k+ \ve j) \not= 0$,
when $|\ve|<1$. Thus,  the inequality $C_k>0$ yields $(C_k+ \ve
j)\geq 0$ for $|\ve| \leq 1$.
\end{proof}
\section{Weyl functions, direct and inverse \\ problem:
the case of the pseudoexponential potentials} \label{direct}
 \setcounter{equation}{0}
Following definitions of the Weyl functions for Sturm-Liouville,
Dirac type and canonical systems on the semiaxis (see, for
instance, \cite{LeS, SaL3} and references therein), we can define
also Weyl functions for system (\ref{0.1}). Namely, let matrices
$C_k>0$ satisfy (\ref{0.2}). Then, a $p \times p$ matrix function
$\vp$ holomorphic in  the lower halfplane $\BC_-$ is said to be a
Weyl function for system (\ref{0.1}) on the semiaxis $k \geq 0$,
 if the  inequality
\begin{equation} \label{1.9}
\sum_{k=0}^\infty [\vp(\lambda)^* \quad
I_p]q(\lambda)^kW_k(\lambda)^*C_kW_k(\lambda)\left[\begin{array}{c}
 \vp(\lambda) \\ I_p
\end{array}
\right]<\infty
\end{equation}
holds, where $q(\lambda)=|\lambda^2|(|\lambda^2|+1)^{-1}$.
\begin{Rk} Similar to the continuous case we have a summation formula:
\begin{equation} \label{1.10}
\sum_{k=0}^r
q(\lambda)^kW_k(\lambda)^*C_kW_k(\lambda)=\frac{|\lambda^2|+1}{i(
\lambda - \ov
\lambda)}\Big(q(\lambda)^{r+1}W_{r+1}(\lambda)^*jW_{r+1}(\lambda)-j\Big).
\end{equation}
Indeed, according to (\ref{0.1}) and (\ref{0.2}) we have
\[
W_{k+1}(\lambda)^*jW_{k+1}(\lambda)=W_{k}(\lambda)^*\Big( I_m +
\frac{i}{\ov \lambda} C_k j \Big)j\Big( I_m -\frac{i}{\lambda} j
C_k \Big)W_{k}(\lambda)
\]
\[
=q(\lambda)^{-1}W_k(\lambda)^*jW_k(\lambda)+\frac{i( \lambda -\ov
\lambda)}{|\lambda^2|}W_k(\lambda)^*C_kW_k(\lambda),
\]
i.e.,
\[
\frac{|\lambda^2|+1}{i(\lambda - \ov \lambda)}\Big(
 q(\lambda)^{k+1}W_{k+1}(\lambda)^*jW_{k+1}(\lambda)-q(\lambda)^{k}W_{k}(\lambda)^*jW_{k}(\lambda)\Big)
\]
\begin{equation} \label{1.11}
=q(\lambda)^kW_k(\lambda)^*C_kW_k(\lambda).
\end{equation}
Formula (\ref{1.11}) yields  (\ref{1.10}).
\end{Rk}
To construct the Weyl function, partition first matrix $\Pi_0$ and
matrix-function $w_A(0,\lambda)$ into  blocks:
\begin{equation} \label{1.12}
\Pi_0=[\Phi \quad \Psi], \quad w_A(0,\lambda)=\left[
\begin{array}{lr}
a(\lambda) & b(\lambda) \\ c(\lambda) & d(\lambda)
\end{array}
\right].
\end{equation}
Similar to the considerations in \cite{GKS1} it follows from
(\ref{0.14}) that
\begin{equation} \label{1.13}
b(\lambda)d(\lambda)^{-1}=-i\Phi^*S_0^{-1}(A^{\times}-\lambda
I_n)^{-1}\Psi, \quad A^{\times}=A+i\Psi\Psi^*S_0^{-1}.
\end{equation}
To calculate $d(\lambda)^{-1}$ here, we use the fact from the
system theory:
\begin{equation} \label{1.14}
\Big(I_p+C(\lambda I_n-A)^{-1}B\Big)^{-1}=I_p- C\big(\lambda
I_n-(A-BC)\big)^{-1}B.
\end{equation}
\begin{Tm} \label{TmDirect}
Let parameter matrices be fixed, assume $S_0>0$, and define $C_k$
by (\ref{0.13}). Then system (\ref{0.1}) is well-defined on the
semiaxis and its unique Weyl function, which satisfies
(\ref{1.9}), takes the form
\begin{equation} \label{1.15}
\vp(\lambda)=-i\Phi^*S_0^{-1}(A^{\times}-\lambda I_n)^{-1}\Psi,
\quad A^{\times}=A+i\Psi\Psi^*S_0^{-1}.
\end{equation}
\end{Tm}
\begin{proof}. By Proposition \ref{Pn3.1} system (\ref{0.1}) is well-defined.
Now, relations (\ref{1.13}) imply $\vp=bd^{-1}$ for the matrix
function $\vp$ given by (\ref{1.15}). According to (\ref{0.22})
we have
\[
w_A(0,\lambda)^*jw_A(0,\lambda) \leq j \quad (\lambda \in \BC_-),
\]
and it follows, in particular, that $d(\lambda)^*d(\lambda) \geq
I_p + b(\lambda)^* b(\lambda)$. Therefore, we get
\begin{equation} \label{1.16}
\vp(\lambda)^* \vp(\lambda)<I_p \quad (\lambda \in \BC_-),
\end{equation}
and so $\vp$ is holomorphic in $\BC_-$. Notice that the equality
$\vp=bd^{-1}$ is equivalent to the formula
\begin{equation} \label{1.17}
\left[\begin{array}{c} \vp(\lambda) \\ I_p
\end{array}
\right]=w_A(0,\lambda)\left[\begin{array}{c} 0 \\ I_p
\end{array}
\right]d(\lambda) ^{-1}.
\end{equation}
Taking into account  (\ref{1.17}) and
$w_A(r+1,\lambda)^*jw_A(r+1,\lambda) \leq j$, we derive from
representation (\ref{0.15}) of $W_{r+1}( \lambda)$ that
\[
[ \vp(\lambda)^* \quad I_p
]W_{r+1}(\lambda)^*jW_{r+1}(\lambda)\left[\begin{array}{c}
 \vp(\lambda) \\ I_p
\end{array}
\right]=|\lambda+i|^{2r+2}|\lambda|^{-2r-2}\big(d(\lambda)^*\big)^{-1}
\]
\begin{equation} \label{1.18}
\times [0 \quad I_p]
w_A(r+1,\lambda)^*jw_A(r+1,\lambda)\left[\begin{array}{c} 0 \\ I_p
\end{array}
\right]d(\lambda) ^{-1}<0.
\end{equation}
By (\ref{1.10}) and (\ref{1.18}) the  inequality
\begin{equation} \label{1.19}
[ \vp(\lambda)^* \quad I_p ]\sum_{k=0}^r
q(\lambda)^kW_k(\lambda)^*C_kW_k(\lambda)\left[\begin{array}{c}
 \vp(\lambda) \\ I_p
\end{array}
\right]<\frac{|\lambda^2|+1}{i( \lambda - \ov \lambda)}I_p
\end{equation}
is true. From (\ref{1.19})  inequality (\ref{1.9}) is immediate,
i.e., $\vp$ defined by  (\ref{1.15}) is a Weyl function.

Let us show that $\vp$ is a unique Weyl function. First notice
that by Proposition \ref{Pn3.2} we have inequality
$W_s^*C_sW_s\geq W_s^*jW_s$. Now, use relation (\ref{1.11}) to
derive inequality $q^sW_s^*jW_s \geq  q^{s-1}W_{s-1}^*jW_{s-1}$.
From the  inequalities above we get
\begin{equation} \label{1.20}
q(\lambda)^kW_k(\lambda)^*C_kW_k(\lambda)\geq j.
\end{equation}
Therefore, the following equality is immediate for any $f \in
\BC^p$:
\begin{equation} \label{1.21}
\sum_{k=0}^\infty f^*[I_p \quad 0
]q(\lambda)^kW_k(\lambda)^*C_kW_k(\lambda)\left[\begin{array}{c}
I_p \\ 0
\end{array}
\right]f =\infty .
\end{equation}
According to (\ref{1.9}) and (\ref{1.21}), the dimension of the
subspace $L \in \BC^m$, such that for all $h \in L$ we have
\begin{equation} \label{1.22}
\sum_{k=0}^\infty h^*
q(\lambda)^kW_k(\lambda)^*C_kW_k(\lambda)h<\infty ,
\end{equation}
equals $p$. Now, suppose  that there is a Weyl function $\wt \vp
\not= \vp$, where $\vp$ is given by (\ref{1.15}). Then the columns
of $\left[\begin{array}{c}
 \vp(\lambda) \\ I_p
\end{array}
\right]$ and the columns of $\left[\begin{array}{c}
 \wt \vp(\lambda) \\ I_p
\end{array}
\right]$ belong to $L$. Therefore, $\dim \, L>p$ for those
$\lambda$, where  $\wt \vp(\lambda)  \not= \vp(\lambda)$, and we
come to a contradiction.
\end{proof}
\begin{Rk} \label{Id}
If $S_0>0$, then by (\ref{0.10})-(\ref{0.13}) we can substitute
parameter matrices $A$, $S_0$ and $\Pi_0$ by the parameter
matrices $S_0^{-\frac{1}{2}}AS_0^{\frac{1}{2}}$, $I_n$ and
$S_0^{-\frac{1}{2}}\Pi_0$, which determine the same system. For
$S_0=I_p$ formula  (\ref{1.15}) takes the form
\begin{equation} \label{1.15'}
\vp(\lambda)=-i\Phi^*(A^{\times}-\lambda I_n)^{-1}\Psi, \quad
A^{\times}=A+i\Psi\Psi^*,
\end{equation}
and we have also $A^{\times}-(A^{\times})^*=i( \Phi  \Phi^*+ \Psi
\Psi^*)$, $\det(A^{\times}-i\Psi\Psi^*)\not=0$.
\end{Rk}
\begin{Ee} \label{pseudo}
Consider the simplest example: $p=1$, $n=1$, $A=a \in \BR$
($a\not=0$), $S_0=1$. From  (\ref{0.9}) and (\ref{0.10}) it
follows that $|\Phi|=|\Psi|$ and
\begin{equation} \label{1.30}
\Pi_k=[\Big(\frac{a+i }{a}\Big)^k \Phi \quad \Big(\frac{a-i
}{a}\Big)^k \Psi ], \quad \Pi_k \Pi_k^*=2
|\Phi|^2\Big(\frac{a^2+1 }{a^2} \Big)^k.
\end{equation}
Now, in view of $S_0=1$, (\ref{0.11}) and the second relation in
(\ref{1.30}) one can check that
\begin{equation} \label{1.31}
S_k=(k \zeta+1)\Big(\frac{a^2+1 }{a^2} \Big)^k, \quad
\zeta=\frac{2|\Phi|^2 }{a^2+1}.
\end{equation}
Finally, using  (\ref{0.13}), (\ref{1.30}) and (\ref{1.31}) we get
the entries $(C_k)_{ij}$ of $C_k$:
\begin{equation} \label{1.32}
(C_k)_{11}=(C_k)_{22}=1+ \zeta |\Phi|^2 (k \zeta+1)^{-1}((k
+1)\zeta+1)^{-1},
\end{equation}
\begin{equation} \label{1.33}
(C_k)_{21}=\ov{(C_k)_{12}}=\Phi \ov \Psi\Big((k
\zeta+1)^{-1}\Big(\frac{a+i }{a-i} \Big)^k- ((k
+1)\zeta+1)^{-1}\Big(\frac{a+i }{a-i} \Big)^{k+1} \Big).
\end{equation}
The Weyl function of system (\ref{0.1}), where the matrices $C_k$
are given by (\ref{1.32}) and (\ref{1.33}), is easily calculated
using (\ref{1.15'}):
\begin{equation} \label{1.34}
\vp(\lambda)=i \ov \Phi  \Psi(\lambda - a-i |\Psi|^2 )^{-1}.
\end{equation}
\end{Ee}

Notice that our matrices $C_k$ belong to the class of the so
called pseudoexponential potentials. An important subclass of the
strictly pseudoexponential potentials, that is, a subclass with an
additional requirement $\s(A) \subset \BC_-$  ($\s$ - spectrum),
have been treated for $p=1$ in \cite{AG1, AG2}. In particular, for
the strictly pseudoexponential subcase the inequality
 $|\vp(\lambda)|<1$ for $\lambda \in \ov{\BC_-}$ is true. On the other hand,  in the simple example above
we have $\s(A)=a \in \BR$ and $|\vp|=1$ for $\lambda=a$.

According to (\ref{1.15}) and (\ref{1.16})  Weyl function $\vp$ is
a rational, strictly proper  and  contractive in $\BC_-$ matrix
function. By the proof of Theorem 9.4 \cite{GKS6} such matrix
functions admit representation (realisation):
\begin{equation} \label{1.23}
\vp(\lambda)=- i \wt \Phi^*(\t-\lambda I_n)^{-1}\wt \Psi,
\end{equation}
where
\begin{equation} \label{1.24}
\t - \t^*=i(\wt \Phi \wt \Phi^*+\wt \Psi \wt \Psi^*).
\end{equation}
Direct calculation shows also that formulas (\ref{1.23}) and
(\ref{1.24}) yield $I_p-\vp^*\vp \geq 0$ for $\lambda \in \BC_-$.
So, realization (\ref{1.23}), (\ref{1.24}) is equivalent to
function being rational, strictly proper  and  contractive in
$\BC_-$.
\begin{Tm}\label{invpr}   Matrix function $\vp$ is the  Weyl function
of some system (\ref{0.1}) determined by the parameter matrices
with $S_0>0$ if and only if it admits representation (\ref{1.23}),
(\ref{1.24}) such that $\det ( \t -i \wt \Psi \wt \Psi^*)\not= 0$.
Then $\vp$ is the Weyl function of some system (\ref{0.1}), where
$C_k>0$. To recover such  system  put
\begin{equation} \label{1.25}
S_0=I_n, \quad A=\t- i\wt \Psi \wt \Psi^*, \quad  \Phi =\wt \Phi ,
\quad \Psi = \wt \Psi, \quad \Pi_0=[\Phi \quad \Psi],
\end{equation}
and  define matrices $C_k$ by formula (\ref{0.13}), where matrices
$\Pi_k$ and $S_k$ ($k>0$) are given by formulas (\ref{0.10}) and
(\ref{0.11}).
\end{Tm}
\begin{proof}. The necessity of theorem's conditions follows from Remark \ref{Id}.
Now, suppose these conditions are fulfilled. Then, from
(\ref{1.24}) and (\ref{1.25}) it follows that the identity
(\ref{0.9}) holds for the parameter matrices. Therefore, system
(\ref{0.1})  is defined. So, by Theorem \ref{TmDirect}
 $\vp$ is the  Weyl function of this system.
\end{proof}

\begin{Rk}  The Weyl functions in the upper halfplane can be treated in a quite similar
way. That is, we define Weyl functions in $\BC_+$ by the
inequality
\begin{equation} \label{vp1.9}
\sum_{k=0}^\infty [I_p \quad \vp(\lambda)^*
]q(\lambda)^kW_k(\lambda)^*C_kW_k(\lambda)\left[\begin{array}{c}
I_p \\ \vp(\lambda)
\end{array}
\right]<\infty .
\end{equation}
Then the Weyl function of system (\ref{0.1}), where matrices $C_k$
are given by (\ref{0.13}) and $S_0>0$, takes the form
\begin{equation} \label{n1.15}
\vp(\lambda)=c(\lambda)a(\lambda)^{-1}=
i\Psi^*S_0^{-1}(A^{\times}-\lambda I_n)^{-1}\Phi, \quad
A^{\times}=A-i\Phi\Phi^*S_0^{-1}.
\end{equation}
\end{Rk}
A  definition of a Weyl function in $\BC_-$ can be also given in a
more general form.
\begin{Dn} \label{WT} Let matrices $C_k>0$ satisfy (\ref{0.2}). Then, a $p \times p$ matrix function $\vp$
holomorphic in $\BC_-$ is said to be a Weyl function for system
(\ref{0.1}) on the semiaxis $k \geq 0$,
 if the following inequality holds:
\begin{equation} \label{2.76}
\sum_{k=0}^\infty [i\vp(\lambda)^* \quad I_p ]q(\lambda)^k K
W_k(\lambda)^*C_kW_k(\lambda)K^*\left[\begin{array}{c}
 -i\vp(\lambda) \\ I_p
\end{array}
\right]<\infty .
\end{equation}
Here $K^*=K^{-1}$ and
$q(\lambda)=|\lambda^2|(|\lambda^2|+1)^{-1}$.
\end{Dn}
When $K$ in the inequality  (\ref{2.76}) equals $I_n$, this
inequality coincides with the inequality  (\ref{1.9}). In general,
the choice of the matrix $K$  is related to the choice of the
domain of the operator corresponding to the Dirac system, and
usually $K$ is chosen so that the Weyl functions are Herglotz
functions. Further we assume that
\begin{equation} \label{1.26}
K=\frac{1}{\sqrt{2}}\left[\begin{array}{lr} I_p & -I_p \\ I_p &
I_p
\end{array}
\right],
\end{equation}
Simple transformations show that the Weyl function $\vp_I$ defined
via (\ref{1.9}) and the Weyl function $\vp_K$ defined via
(\ref{2.76}) and  (\ref{1.26}) are connected by the relation
\begin{equation} \label{1.27}
\vp_K=-i(I_p- \vp_I)(I_p+ \vp_I)^{-1}.
\end{equation}
From (\ref{1.27}) it follows that
\begin{equation} \label{1.28}
\vp_K(\lambda)-\vp_K(\lambda)^*=-2i(I_p+
\vp_I(\lambda)^*)^{-1}(I_p-\vp_I(\lambda)^*\vp_I(\lambda))(I_p+
\vp_I(\lambda))^{-1},
\end{equation}
$\lambda \in \BC_-$. Thus, according to (\ref{1.16}) and
(\ref{1.28}) $\vp_K$ is a Herglotz function with a non-positive
imaginary part in $\BC_-$.

\section{Weyl functions, direct and inverse problem on the interval:
general case} \label{gencase}
 \setcounter{equation}{0}
In this section we shall consider the self-adjoint matrix discrete
Dirac type system (\ref{0.1}) on the interval $0 \leq k \leq N$.
We assume that (\ref{0.2}) holds and $C_k>0$. It was  shown in
Proposition \ref{equiv}  that $C_k>0$ yields
 $C_k=(U_k^*)^{-1}R_k^2 U_k^{-1}$, where $R_k=R_k^*$ and $R_kjR_k=j$. Hence, we get
\[
C_k=(U_k^*)^{-1}(R_k^2+R_kjR_k)
U_k^{-1}-j=(U_k^*)^{-1}R_k(I_m+j)R_k U_k^{-1}-j
\]
\[
=2 \wh \b(k)^* \wh \b(k)-j, \quad 0 \leq k \leq N,
\]
where
\begin{equation} \label{RU}
\wh \b(k)=[I_p\quad 0]R_k U_k^{-1}, \quad \wh \b(k) j \wh
\b(k)^*=I_p.
\end{equation}
Further we shall use these relations:
\begin{equation} \label{2.2}
C_k=2 \wh \b(k)^* \wh \b(k)-j, \quad \wh \b(k) j \wh \b(k)^*=I_p,
\quad 0 \leq k \leq N.
\end{equation}
\begin{Rk} \label{proC}
Relations (\ref{2.2}) are equivalent to the relations $C_k>0$ and
(\ref{0.2}). Indeed, we have just derived (\ref{2.2}) from $C_k>0$
and $C_kjC_k=j$, and vice versa: direct calculation shows that
(\ref{2.2}) yields (\ref{0.2}). To derive also from (\ref{2.2})
the inequality $C_k>0$, choose a matrix $\wt \b(k)$ such that
\begin{equation}\label{2.14'}
\wt \b(k) j \wh \b(k)^*=0, \quad \wt \b(k) j \wt \b(k)^*=-I_p.
\end{equation}
Notice, that in view of the second relations in (\ref{2.2}) the
 maximal subspace, which is $j$-orthogonal to the rows of $\wh \b_k$,
proves to be $p$-dimensional and $j$-negative, i.e., $\wt \b(k)$
always exists. According to (\ref{2.2}) and (\ref{2.14'}) we have
\begin{equation}\label{2.14n}
 \left[
\begin{array}{c} \wh \b(k)
\\ \wt \b(k)
\end{array}
\right]j \left[
\begin{array}{c} \wh \b(k)
\\ \wt \b(k)
\end{array}
\right]^*=j=\left[
\begin{array}{c} \wh \b(k)
\\ \wt \b(k)
\end{array}
\right]^*j
 \left[
\begin{array}{c} \wh \b(k)
\\ \wt \b(k)
\end{array}
\right].
\end{equation}
Finally, by the first relations in (\ref{2.2}) and by
(\ref{2.14n}) we obtain
\begin{equation}\label{2.14''}
C_k=\wh \b(k)^*  \wh \b(k)+ \wt \b(k)^*  \wt \b(k)=\left[
\begin{array}{c} \wh \b(k)
\\ \wt \b(k)
\end{array}
\right]^* \left[
\begin{array}{c} \wh \b(k)
\\ \wt \b(k)
\end{array}
\right]>0.
\end{equation}
\end{Rk}
From the second relation in (\ref{2.2}) for $k\geq 0$ it follows
also that
\begin{equation} \label{2.2''}
\det\Big( \wh \b(k) j \wh \b(k+1)^*\Big)\not=0, \quad  0 \leq k
\leq N-1 .
\end{equation}
Indeed, if (\ref{2.2''}) does not hold, we have $\wh \b(k) j \wh
\b(k+1)^*f=0$ for some $f\not= 0$. Then, in view of the second
relations in (\ref{2.2}) for $k\geq 0$, we see that the linear
span of the rows of $\b_k$ and of $f^* \b_{k+1}$ forms a
$p+1$-dimensional $j$-positive subspace of $\BC_m$, which is
impossible.

Notice also that according to  (\ref{0.6})  and (\ref{0.8}) we
have
\begin{equation} \label{RU1}
U_{k}^{-1}=-i j R_{k-1}U_{k-1}^{-1}, \quad U_k^*=-i
R_{k-1}U_{k-1}^{-1}j \quad (k>0).
\end{equation}
Compare (\ref{RU})  and (\ref{2.2}) to see that a matrix $\wh
\b(k-1)$,  which satisfies the equalities
\[
C_{k-1}=2 \wh \b(k-1)^* \wh \b(k-1)-j, \quad \wh \b(k-1) j \wh
\b(k-1)^*=I_p,
\]
coincides with the upper block row of $-i  R_{k-1}U_{k-1}^{-1}$ up
to a $p \times p$ unitary factor $\wh u(k-1)$. Moreover, $\wt
\b(k-1)$ given by  (\ref{2.14'}) defines the lower block row of
$-i  R_{k-1}U_{k-1}^{-1}$ up to some $p \times p$ unitary factor
$\wt u(k-1)$. Taking into account the second equalities in
(\ref{n.2}) and (\ref{RU1}), we get
\begin{equation} \label{RU2}
R_k=\left( D(k-1)\left[
\begin{array}{c} \wh \b(k-1)
\\ \wt \b(k-1)
\end{array}
\right]j C_k j \left[
\begin{array}{cc} \wh \b(k-1)^*
& \wt \b(k-1)^*
\end{array}
\right] D(k-1)^*\right)^{\frac{1}{2}},
\end{equation}
where
\begin{equation} \label{RU3}
D(k-1)= \left[
\begin{array}{cc} \wh u(k-1) & 0 \\
0 & \wt u(k-1)
\end{array}
\right].
\end{equation}
\begin{Rk} \label{RK}
We see that the two matrices $C_k$ and $C_{k-1}$ determine $R_k$
by (\ref{RU2}) up to a unitary block diagonal matrix $D(k-1)$.
\end{Rk}

Similar to the continuous case,  Weyl functions  of the discrete
system on the interval are defined via M\"obius
(linear-fractional) transformation
\begin{equation}\label{2.5}
\varphi (\lambda ) = i \bigl( {\cal W}_{21}(\lambda )R(\lambda )+
{\cal W}_{22}(\lambda )Q(\lambda )\bigr) \bigl( {\cal
W}_{11}(\lambda ) R(\lambda )+{\cal W}_{12}(\lambda )Q(\lambda
)\bigr)^{-1},
\end{equation}
where $R$ and $Q$ are $p \times p$ analytical functions in the
neighbourhood of $\lambda=-i$, and
\begin{equation} \label{2.4}
{\cal W}(\lambda)= \{{\cal W}_{ij}(\lambda) \}_{i,j=1}^2=K
W_{N+1}( \ov{\lambda})^*,
\end{equation}
Here,  coefficients ${\cal W}_{ij}$ of the M\"obius
transformation are the $p \times p$ blocks of ${\cal W}$,  the
matrix $K$ is given by (\ref{1.26}) and
\begin{equation}\label{2.5'}
K^*=K^{-1}, \quad KjK^*=J, \quad J=\left[
\begin{array}{cc}
0 & I_{p}  \\ I_{p} & 0
\end{array}
\right].
\end{equation}
It would be convenient to put $\b(k):=\wh \b(k) K^*$ and rewrite
(\ref{2.2}) as
\begin{equation} \label{2.4''}
 C_k=2 K^* \b(k)^*\b(k)K-j, \quad \b(k) J \b(k)^*=I_p, \quad  0 \leq k \leq N.
\end{equation}
We shall need the following analog (for the self-adjoint case) of
Theorem 3.4 \cite{SaA8}.
\begin{Tm} \label{Tm2.2} Suppose
$W$ $\, (W_0(\lambda)=I_m)$ is the fundamental solution of  system
(\ref{0.1}), which   satisfies conditions (\ref{2.4''}). Suppose
also that a $p \times p$ matrix function $\vp$ is given by
formulas (\ref{2.5}) and (\ref{2.4}), where
\begin{equation} \label{2.4'}
\det \bigl( {\cal W}_{11}(-i ) R(-i )+{\cal W}_{12}(-i )Q(-i
)\bigr)\not= 0.
\end{equation}
 Then system (\ref{0.1}) satisfies  (\ref{0.2}), $C_k>0$ ($k \geq 0$), and
the inequalities
\begin{equation} \label{2.4n}
\det\Big( \b(k) J  \b(k+1)^*\Big)\not=0, \quad  0 \leq k \leq N-1
\end{equation}
hold. Moreover, system (\ref{0.1}) is uniquely recovered from the
first $N+1$ Taylor coefficients $\{ \a_k \}_{k=0}^N$ of $\,i
\displaystyle{\vp \left(i \Big( \frac{z+1}{z-1} \Big) \right) }$
at $z=0$ by the following procedure.

First, introduce $(N+1)p \times p$ matrices $\Phi_1$, $\Phi_2:$
\begin{equation}\label{2.11}
\Phi_1 = \left[
\begin{array}{c}
I_{p}  \\ I_p \\ \cdots \\ I_{p}
\end{array}
\right], \quad \Phi_2 =  \left[
\begin{array}{l}
\a_0  \\ \a_0+ \a_1 \\ \cdots \\ \a_0+ \a_1 + \ldots + \a_N
\end{array}
\right].
\end{equation}
Then, introduce an $(N+1)p \times 2p$ matrix $\Pi$ and an $(N+1)p
\times (N+1)p$ block lower triangular matrix $A$ by the blocks $:$
$\quad \Pi=[\Phi_1 \quad \Phi_2]$,
\begin{equation}\label{2.12}
\begin{array}{lcllcl}
 A:=A(N) &=& \left\{ a_{j-k}^{\,} \right\}_{k,j=0}^N,
           & a_r  &=&  \left\{ \begin{array}{lll}
                                  0 \, & \mbox{ for }& r > 0   \\
                                 \displaystyle{\frac{i}{ {\, 2 \,}}} \,
                                 I_p
                                   & \mbox{ for }& r = 0   \\
                                 \, i \, I_p
                                   & \mbox{ for }& r < 0
                           \end{array} \right. \end{array}.
\end{equation}
Next, we recover $(N+1)p \times (N+1)p$ matrix $S$ as a unique
solution of the matrix identity
\begin{equation} \label{2.13}
AS-SA^*=i \Pi J \Pi^*.
\end{equation}
This solution is invertible and positive, i.e., $S>0$. Finally,
matrices $\b(k)^*\b(k)$ are easily recovered from the formula
\begin{equation}\label{2.14}
\Pi^*S^{-1}\Pi=B^*B, \quad B: =B(N)= \left[
\begin{array}{c}
\b(0) \\ \b(1) \\ \cdots \\ \b(N)
\end{array}
\right].
\end{equation}
Now, matrices $C_k$ and system (\ref{0.1}) are defined via the
first equality in (\ref{2.4''}).
\end{Tm}
\begin{proof}. Step 1.  According to Remark \ref{proC} relations $C_k>0$
and (\ref{0.2}) follow from (\ref{2.4''}).  Relations (\ref{2.4n})
follow from (\ref{2.2''}). Now, let
\begin{equation}\label{2.15}
K(r) = \left[
\begin{array}{c}
K_0(r) \\ K_1(r) \\ \cdots \\ K_r(r)
\end{array}
\right],
\end{equation}
where $K_l(r)$ are $p \times (r+1)p$ matrices of the form
\begin{equation}\label{2.16}
K_l(r)=i \b(l)J[\b(0)^* \ldots \b(l-1)^* \quad \b(l)^*/2 \quad 0
\ldots 0].
\end{equation}
From (\ref{2.14})-(\ref{2.16}) it follows that
\begin{equation}\label{2.17}
K(r)-K(r)^*=i B(r) J B(r)^*.
\end{equation}
By induction we shall show in the next step that $K$ is similar to
$A$:
\begin{equation}\label{2.18}
K(r)=V_-(r)A(r)V_-(r)^{-1} \quad (0 \leq r \leq N),
\end{equation}
where $V_-(r)^{\pm 1}$ are block lower triangular matrices. Taking
into account (\ref{2.18}) and multiplying both sides  of
(\ref{2.17}) by $V_-(r)^{-1}$ from the left  and by $\big(V_-(r)^*
\big)^{-1}$ from the right,  we get
\begin{equation}\label{2.19}
A(r)S(r)-S(r)A(r)^*=i \Pi(r) J \Pi(r)^*,
\end{equation}
\begin{equation}\label{2.20}
S(r):=V_-(r)^{-1}\big(V_-(r)^* \big)^{-1}, \quad
\Pi(r):=V_-(r)^{-1}B(r).
\end{equation}
Moreover,  Step 3 will show that matrix $V_-(N)$ can be chosen so
that the equality
\begin{equation}\label{2.21}
\Pi=[\Phi_1 \quad \Phi_2]=V_-(N)^{-1}B(N)
\end{equation}
holds, i.e., $\Pi= \Pi(N)$. (Here $\Phi_1$ and $\Phi_2$ are given
by (\ref{2.11}).)

Identities (\ref{2.19}) have unique solutions $S(r)$ as the
spectra of $A(r)$ and $A(r)^*$ do not intersect. (The statement
follows from the rewriting of (\ref{2.19}) in the form
\[
S(r)(A(r)^*-\lambda I)^{-1}-(A(r)-\lambda I)^{-1}S(r)
\]
\[
=
i(A(r)-\lambda I)^{-1}\Pi(r) J \Pi(r)^* (A(r)^*-\lambda I)^{-1},
\]
and from the following integration of  both sides of the obtained
identity  along a contour, such that the spectra of $A$ is inside
and the spectra of $A^*$ outside it.  In particular, by
(\ref{2.13}) and (\ref{2.19}) one can see that $S:=S(N)$. Hence,
we derive from (\ref{2.20}) and (\ref{2.21}) that $S>0$ and the
first equality in (\ref{2.14}) holds. It remains only to prove
(\ref{2.18}) and (\ref{2.21}).

Step 2. Now, we shall consider block lower triangular matrices
$V_-(k)$ $\, (0 \leq k \leq N)$:
\begin{equation}\label{2.22}
V_-(0)=v_-(0)=\b_1(0), \quad V_-(k)= \left[
\begin{array}{cc}
V_-(k-1) & 0 \\ X(k) & v_-(k)
\end{array}
\right] \quad (k>0),
\end{equation}
where $v_-(k)$ are $p \times p$ matrices, where $\b_1(k)$ and
$\b_2(k)$ are $p \times p$ blocks of $\b(k)=[\b_1(k) \quad
\b_2(k)]$, and where $X(k)=[X_0(k) \quad \wt X(k) ]$ are $p \times
k p$ matrices. Here $X_0(k)$ are  arbitrary $p \times p$ blocks,
and the matrices $\wt X(k)$, $v_-(k)$ are given by the formulas
\[
\wt X(k)=i \Big( \b(k)J[ \b(0)^* \ldots \b(k-1)^*]V_-(k-1)\left[
\begin{array}{c}
I_{(k-1)p}\\ 0
\end{array}
\right]-v_-(k)[I_p \ldots I_p] \Big)
\]
\begin{equation}\label{2.23}
\times \Big(A(k-2)+ \frac{i}{2} I_{(k-1)p}\Big)^{-1}, \quad
v_-(k)=\b(k)J\b(k-1)^* v_-(k-1).
\end{equation}
According to (\ref{2.12}) we have $A(0)=(i/2)I_p$. From the second
relation in (\ref{2.4''}) and definitions (\ref{2.15}) and
(\ref{2.16}) it is immediate that $K(0)=(i/2)I_p$, and so
(\ref{2.18}) is valid for $r=0$. Assume that (\ref{2.18}) is true
for $r=k-1$, and let us show that (\ref{2.18}) is true for $r=k$
too. It is easy to see that
\begin{equation}\label{2.25}
 V_-(k)^{-1}= \left[
\begin{array}{cc}
V_-(k-1)^{-1} & 0 \\ -v_-(k)^{-1} X(k)V_-(k-1)^{-1}  & v_-(k)^{-1}
\end{array}
\right].
\end{equation}
Then, in view of definitions (\ref{2.12}) and (\ref{2.22}), our
assumption implies
\begin{equation}\label{2.26}
V_-(k)A(k)V_-(k)^{-1}=\left[
\begin{array}{cc}
K(k-1) & 0 \\ Y(k) & \frac{i}{2} I_{p}
\end{array}
\right],
\end{equation}
where $Y(k)=\Big[\big( X(k)A(k-1)+iv_-(k)[I_p \ldots I_p]\big)
\quad \, \frac{i}{2}v_-(k) \Big]$
\[
\times \left[
\begin{array}{c}
V_-(k-1)^{-1}  \\ -v_-(k)^{-1} X(k)V_-(k-1)^{-1}
\end{array}
\right].
\]
Rewrite the product on the right-hand side of the last formula as
\begin{equation}\label{2.27}
Y(k)=\Big( X(k)\big(A(k-1)- \frac{i}{2}I_{kp} \big)+iv_-(k)[I_p
\ldots I_p]\Big)V_-(k-1)^{-1}.
\end{equation}
From (\ref{2.12}) and (\ref{2.27}) it follows that
\begin{equation}\label{2.28}
Y(k)=\Big[ \Big(\wt X(k)\big(A(k-2)+ \frac{i}{2}I_{(k-1)p}
\big)+iv_-(k)[I_p \ldots I_p] \Big) \quad
iv_-(k)\Big]V_-(k-1)^{-1}.
\end{equation}
Notice that the sequence $[I_p \ldots I_p]$ of identity matrices
in (\ref{2.28}) is one block smaller than in (\ref{2.27}). By
(\ref{2.23}) and (\ref{2.28}) we have
\[
Y(k)=i \b(k)J \Big[ [ \b(0)^* \ldots \b(k-1)^*]V_-(k-1)\left[
\begin{array}{c}
I_{(k-1)p}\\ 0
\end{array}
\right] \quad \, \b(k-1)^* v_-(k-1) \Big]
\]
\begin{equation}\label{2.29}
\times V_-(k-1)^{-1}.
\end{equation}
Finally, formulas (\ref{2.22}) and (\ref{2.29}) imply
\begin{equation}\label{2.30}
Y(k)=i \b(k) J [ \b(0)^* \ldots \b(k-1)^*] \quad (k>0).
\end{equation}
According to the second relation in (\ref{2.4''}) and formulas
(\ref{2.16}) and (\ref{2.30}) we get
\begin{equation}\label{2.31}
\Big[Y(k) \quad  \frac{i}{2}I_{p} \Big]=K_k(k).
\end{equation}
Now, using  (\ref{2.15})  and (\ref{2.31}), one can see that the
right-hand side of (\ref{2.26}) equals $K(k)$. Thus, (\ref{2.18})
is true for $r=k$ and, therefore, it is true for all $0 \leq r
\leq n$.

Step 3. To derive (\ref{2.21}) we shall first prove that the
matrices $V_-(r)$ given by (\ref{2.22}) and (\ref{2.23}) can be
chosen so that
\begin{equation}\label{2.32}
 V_-(r)^{-1}B_1(r)=\left[
\begin{array}{c}
I_{p}   \\ \cdots \\ I_{p}
\end{array}
\right], \quad B_1(r):= B(r)\left[
\begin{array}{c}
 I_p \\0
\end{array}
\right]=\left[
\begin{array}{c}
 \b_1(0) \\ \cdots \\   \b_1(r)
\end{array}
\right].
\end{equation}
In other words, the   blocks $X_0(r)$, arbitrary till now, can be
chosen so. Indeed, by the definition in (\ref{2.14}) and the first
equality in (\ref{2.22}) formula (\ref{2.32}) is true for $r=0$.
Assume that (\ref{2.32}) is true for $r=k-1$. Then, from
(\ref{2.25}) it follows that (\ref{2.32}) is true for $r=k$, if
only
\begin{equation}\label{2.33}
-v_-(k)^{-1} X(k)\left[
\begin{array}{c}
I_{p}   \\ \cdots \\ I_{p}
\end{array}
\right]+v_-(k)^{-1}\b_1(k)=I_p.
\end{equation}
It implies that we get equality (\ref{2.32})  for $r=k$ by letting
\begin{equation}\label{2.34}
X_0(k)=\b_1(k)-v_-(k)-\wt X(k)\left[
\begin{array}{c}
I_{p}   \\ \cdots \\ I_{p}
\end{array}
\right].
\end{equation}
Hence, by a proper choice of the matrices $X_0(r)$ we  obtain
(\ref{2.32}) for all $r \leq N$.

It remains to prove that
\begin{equation}\label{2.35}
 V_-(N)^{-1}B_2(N)=\Phi_2, \quad B_2(N):=\left[
\begin{array}{c}
 \b_2(0) \\ \cdots \\   \b_2(N)
\end{array}
\right].
\end{equation}
For that purpose we shall consider the matrix function
$W_{N+1}(\lambda)$, which is used in (\ref{2.4}) to define the
coefficients of the M\"obius transformation (\ref{2.5}).  Namely,
we shall prove the  transfer matrix function representation of
$W_{N+1}(\lambda)$:
\begin{equation}\label{2.36}
W_{N+1}(\lambda)=\left(\frac{\lambda+i}{\lambda}\right)^{N+1}K^*
w_A \Big(N, -\frac{ \lambda}{2}\Big)K,
\end{equation}
where
\begin{equation}\label{2.37}
w_A(r, \lambda)=I_{2p}-i J \Pi(r)^*S(r)^{-1}\big(A(r)- \lambda
I_{(r+1)p} \big)^{-1} \Pi(r).
\end{equation}
Identity (\ref{2.19}) implies a similar to (\ref{0.22}) equality
\[
w_A(r, \mu)^*Jw_A(r,\lambda)=J
\]
\begin{equation}\label{2.38}
+i(\ov \mu - \lambda) \Pi(r)^*\big(A(r)^*-\ov \mu
I_{(r+1)p}\big)^{-1}S(r)^{-1}\big(A(r) - \lambda
I_{(r+1)p}\big)^{-1}\Pi(r).
\end{equation}
Moreover, according to factorization theorem 4 from \cite{SaL1}
(see also \cite{SaL3}, p. 188) we have
\[
w_A(r, \lambda)=\Big(I_{2p} -i J
\Pi(r)^*S(r)^{-1}P^*\big(PA(r)P^*- \lambda I_p
\big)^{-1}\big(PS(r)^{-1}P^*\big)^{-1}
\]
\begin{equation}\label{2.39}
\times P S(r)^{-1}\Pi(r) \Big) w_A(r-1, \lambda), \quad P=[0 \quad
\ldots \quad 0 \quad I_p].
\end{equation}
Taking into account (\ref{2.12}), (\ref{2.20}), and (\ref{2.22})
we obtain
\begin{equation}\label{2.40}
\big(PA(r)P^*- \lambda I_p \big)^{-1}=(\frac{i}{2}- \lambda
)^{-1}I_p, \quad PS(r)^{-1}P^*=v_-(r)^*v_-(r),
\end{equation}
\begin{equation}\label{2.41}
P S(r)^{-1}\Pi(r)=v_-(r)^*P B(r)=v_-(r)^* \b(r).
\end{equation}
Substitute (\ref{2.40}) and (\ref{2.41}) into (\ref{2.39}) to get
\begin{equation}\label{2.42}
w_A(r, \frac{\lambda}{2})=\Big(I_{2p} -\frac{2i}{i- \lambda}J
\b(r)^*\b(r) \Big)w_A(r-1, \frac{\lambda}{2}).
\end{equation}
From the definitions (\ref{2.12}), (\ref{2.20}), and (\ref{2.37})
we also easily derive
\begin{equation}\label{2.43}
w_A(0, \frac{\lambda}{2})=I_{2p} -\frac{2i}{i-
\lambda}JB(0)^*B(0)=I_{2p} -\frac{2i}{i- \lambda}J \b(0)^*\b(0).
\end{equation}
On the other hand  system (\ref{0.1}) with additional conditions
(\ref{2.4''}) can be rewritten as
\begin{equation}\label{2.44}
W(r+1, \lambda)=\frac{\lambda +
i}{\lambda}\Big(I_{2p}-\frac{2i}{i+ \lambda}j K^* \b(r)^*\b(r)K
\Big)W(r, \lambda).
\end{equation}
In view of the normalization $W(0)=I_{2p}$, formulas
(\ref{2.42})-(\ref{2.44}) imply  (\ref{2.36}).

From (\ref{2.36}) and (\ref{2.38}) it follows that
\begin{equation}\label{2.45}
W(N+1, \lambda)jW(N+1, \ov
\lambda)^*=\left(\frac{\lambda+i}{\lambda}\right)^{N+1}\left(\frac{\lambda-i}{\lambda}\right)^{N+1}j.
\end{equation}
Let us include  functions $\vp$ into consideration. Introduce
\begin{equation}\label{2.46}
{\cal A}(\lambda):=\left| \frac{\lambda}{\lambda+i}
\right|^{2N+2}[i \vp(\lambda)^* \quad I_p ] KW(N+1,
\lambda)^*jW(N+1, \lambda)K^*\left[
\begin{array}{c}
-i \vp(\lambda) \\ I_p
\end{array}
\right].
\end{equation}
According to (\ref{2.5}), (\ref{2.4}), and (\ref{2.45}) we have
\[
{\cal A}(\lambda)=\left| \frac{\lambda-i}{\lambda} \right|^{2N+2}
\bigl( \big({\cal W}_{11}(\lambda ) R(\lambda )+{\cal
W}_{12}(\lambda )Q(\lambda )\big)^*\bigr)^{-1}
\]
\begin{equation}\label{2.47}
\times \Big(R(\lambda)^*R(\lambda)-
Q(\lambda)^*Q(\lambda)\Big)\bigl( {\cal W}_{11}(\lambda )
R(\lambda )+{\cal W}_{12}(\lambda )Q(\lambda )\bigr)^{-1}.
\end{equation}
By  (\ref{2.4'}) and (\ref{2.47})  ${\cal A}$ is bounded in the
neighborhood of $\lambda=-i$:
\begin{equation}\label{2.48}
\|{\cal A}(\lambda)\|=O(1) \quad {\mathrm{for}} \quad \lambda \to
-i.
\end{equation}
Now, substitute  (\ref{2.36}) and (\ref{2.38}) into (\ref{2.46})
to obtain
\[
{\cal A}(\lambda)=[ i\vp(\lambda)^* \quad I_p ]
\Big(J+\frac{i}{2}(\lambda - \ov \lambda )
\Pi(N)^*\big(A(N)^*+\frac{\ov \lambda}{2}
I_{(N+1)p}\big)^{-1}S(N)^{-1}
\]
\begin{equation}\label{2.49}
\times \big(A(N) + \frac{\lambda}{2}
I_{(N+1)p}\big)^{-1}\Pi(N)\Big)\left[
\begin{array}{c}
-i\vp(\lambda) \\ I_p
\end{array}
\right].
\end{equation}
Notice that $S(N)>0$. Hence, formulas (\ref{2.48}) and
(\ref{2.49}) imply that
\begin{equation}\label{2.50}
\left\| \big(A(N) +\frac{\lambda}{2}
I_{(N+1)p}\big)^{-1}\Pi(N)\left[
\begin{array}{c}
-i \vp(\lambda) \\ I_p
\end{array}
\right]\right\|=O(1) \quad {\mathrm{for}} \quad \lambda \to -i.
\end{equation}
Recall that $\Pi(N)=V_-(N)^{-1}B(N)$ and that $A(N)$ is denoted by
$A$. Now, represent  $\Pi(N)$ in the block form
\begin{equation}\label{2.50'}
\Pi(N)=[\Phi_1(N) \quad \Phi_2(N)], \quad
\Phi_k(N)=V_-(N)^{-1}B_k(N) \quad (k=1,2).
\end{equation}
According to (\ref{2.11}) and (\ref{2.32}) we have
$\Phi_1(N)=\Phi_1$. Hence, multiplying the matrix function on the
left-hand side of (\ref{2.50}) by $i\Big(\Phi_1^*\big(A +
\frac{\lambda}{2} I_{(N+1)p}\big)^{-1}\Phi_1 \Big)^{-1}\Phi_1^*$
we derive
\[\left\|
\vp(\lambda)+i\Big(\Phi_1^*\big(A +\frac{\lambda}{2}
I_{(N+1)p}\big)^{-1}\Phi_1 \Big)^{-1}\Phi_1^*\big(A +
\frac{\lambda}{2} I_{(N+1)p}\big)^{-1}\Phi_2(N)\right\|
\]
\begin{equation}\label{2.51}
=O\left(\left\|\Big(\Phi_1^*\big(A + \frac{\lambda}{2}
I_{(N+1)p}\big)^{-1}\Phi_1 \Big)^{-1}\right\|\right) \quad
{\mathrm{for}} \quad \lambda \to -i.
\end{equation}
The matrix $A + \frac{\lambda}{2} I_{(N+1)p}$ is easily inverted
explicitly (see, for instance, formula (1.10) in \cite{SaAtepl}).
As a result one obtains
\begin{equation}\label{2.52}
\Phi_1^*\big( A + \frac{\lambda}{2}
I_{(N+1)p}\big)^{-1}=\frac{2}{i+\lambda}[\wh q^N \quad \wh q^{N-1}
\quad \ldots \quad \wh q \quad I_p ], \quad  \wh q:=\frac{\lambda
-i}{\lambda +i}I_p.
\end{equation}
Moreover, we get
\begin{equation}\label{2.53}
\Phi_1^*\big( A + \frac{\lambda}{2}
I_{(N+1)p}\big)^{-1}\Phi_1=\frac{2} {i+\lambda}\Big(\wh
q^{N+1}-I_p \Big)\Big( \wh q-I_p \Big)^{-1}.
\end{equation}
 Let $\displaystyle{\lambda =i \Big( \frac{z+1}{z-1}
\Big)}$, i.e., $\displaystyle{z =\Big( \frac{\lambda
+i}{\lambda-i} \Big)}$. Then, we derive from (\ref{2.53}) that
\begin{equation}\label{2.53'}
\Big(\Phi_1^*\big(A + \frac{\lambda}{2} I_{(N+1)p}\big)^{-1}\Phi_1
\Big)^{-1}=\big(-iz^{N+1}+O(z^{2N+2})\big)I_p \quad (z \to 0).
\end{equation}
Taking into account (\ref{2.52}) and (\ref{2.53'}), we rewrite
(\ref{2.51}) as
\begin{equation}\label{2.54}
\left\| \vp\Big(i \Big( \frac{z+1}{z-1} \Big)\Big)+i(1-z)[I_p
\quad zI_p \quad z^2I_p \quad \ldots]\Phi_2(N)\right\| =O(z^{N+1})
\end{equation}
for $z \to 0$. From (\ref{2.11}) and (\ref{2.54}) it follows that
$\Phi_2(N)=\Phi_2$, i.e., (\ref{2.35}) is true. As
$\Phi_1(N)=\Phi_1$ and $\Phi_2(N)=\Phi_2$, so $\Pi(N)=\Pi$ and
formula (\ref{2.21}) is finally proved.
\end{proof}
\begin{Dn} \label{genWeyl}
Let matrices $C_k$ satisfy (\ref{2.4''}). Then, a $p \times p$
matrix function $\vp$ holomorphic in $\BC_-$ is said to be a Weyl
function for system (\ref{0.1}) on the interval $0 \leq k \leq N$,
 if $\vp$ admits representation (\ref{2.5}), where a pair $R$, $Q$ is meromorphic in $\BC_-$,
well-defined at $\lambda =-i$, and nonsingular with $j$-property,
i.e.,
\begin{equation} \label{2.55}
R(\lambda)^*R(\lambda)+Q(\lambda)^*Q(\lambda)>0, \quad
R(\lambda)^*R(\lambda) \leq Q(\lambda)^*Q(\lambda).
\end{equation}

The set of Weyl functions is denoted by ${\cal N}(N)$.
\end{Dn}
Using notation (\ref{2.4}), we deduce from (\ref{1.10})  the
inequality
\begin{equation}\label{2.55'}
q(\lambda)^{N+1}{\cal W}(\lambda)j{\cal W}(\lambda)^*\leq J, \quad
\lambda \in \BC_-.
\end{equation}
According to \cite{Po} we can change the order of  factors in
(\ref{2.55'}):
\begin{equation}\label{2.56}
q(\lambda)^{N+1}{\cal W}(\lambda)^*J{\cal W}(\lambda) \leq j,
\quad \lambda \in \BC_-.
\end{equation}
Moreover, after excluding $\lambda=-i$ the inequality is strict
\begin{equation}\label{2.56'}
q(\lambda)^{N+1}{\cal W}(\lambda)^*J{\cal W}(\lambda) < j, \quad
\lambda \in \BC_-\backslash-i.
\end{equation}
In view of (\ref{0.1}), (\ref{2.4}), (\ref{2.14n}), and
(\ref{2.14''}) at $\lambda=-i$ we get
\begin{equation}\label{2.56n}
{\cal W}(-i)=KW_{N+1}(i)^*=(-2)^{N+1}K \prod_{k=0}^N (\wt
\b_k^*\wt \b_k j).
\end{equation}
From the second relations in (\ref{2.14'}) and from (\ref{2.55})
we, analogously to the proof of (\ref{2.2''}), derive:
\begin{equation}\label{2.56nn}
\det \, [I_p \quad I_p] j \wt \b_0^* \not=0, \quad \det \, \wt
\b_k j \wt \b_{k+1}^* \not=0, \quad \det \, \wt \b_N j \left[
\begin{array}{c}
R(-i) \\ Q(-i)
\end{array}
\right] \not=0.
\end{equation}
By (\ref{2.56n}) and  (\ref{2.56nn}) the   next proposition is
valid.
\begin{Pn} \label{ineq}
Let the pair $R$, $Q$ satisfy (\ref{2.55}) Then inequality
(\ref{2.4'}) is fulfilled.
\end{Pn}
By  Proposition  \ref{ineq} and the proof of Theorem \ref{Tm2.2}
we get a corollary.
\begin{Cy} \label{CyTR} Weyl functions of system (\ref{0.1}),
which satisfies conditions (\ref{2.4''}), are Herglotz functions
and admit the Taylor representation
\begin{equation}\label{2.57'}
 \vp\Big(i \Big( \frac{z+1}{z-1}
\Big)\Big)=-i\Big(\psi_0+(\psi_1-\psi_0)z+\ldots+(\psi_{N}-\psi_{N-1})z^N\Big)+O(z^{N+1})
,
\end{equation}
where $z \to 0$ and the $p \times p$ matrices $\psi_k$ are the
blocks of
\begin{equation}\label{2.57''}
\Phi_2=\left[
\begin{array}{c}
 \psi_0 \\ \psi_1 \\ \ldots \\ \psi_N
\end{array}
\right]=V_-(N)^{-1}B(N).
\end{equation}
\end{Cy}
\begin{proof}.
From (\ref{2.5}), (\ref{2.55}) and (\ref{2.56}) it follows that
\[
[I_p \quad i\vp^*]J  \left[
\begin{array}{c}
 I_p \\ -i \vp
\end{array}
\right]\leq 0,
\]
i.e., $\Im \vp(\lambda) \leq 0$ for $\lambda \in \BC_-$, and so
$\vp$ is a Herglotz function.

By  Proposition  \ref{ineq} the Weyl functions satisfy conditions
of  Theorem \ref{Tm2.2}. Then, by the second relation in
(\ref{2.11}) we have representation (\ref{2.57'}) of $\vp$ via the
blocks of $\Phi_2$. By the proof of Theorem \ref{Tm2.2} we get
also $\Phi_2=\Phi_2(N)$, i.e., (\ref{2.57''}) holds. Here $V_-(N)$
and $B(N)$ are recovered from the matrices $\b(k)$ and do not
depend on $\vp$.
\end{proof}
\begin{Rk} \label{InvPr}
As Weyl functions $\vp$ satisfy conditions of  Theorem
\ref{Tm2.2}, so the procedure given in Theorem \ref{Tm2.2}
provides a recovery of system (\ref{0.1}) from a Weyl function
(i.e., provides a solution of the inverse problem).
\end{Rk}
The following proposition is also true
\begin{Pn} \label{imbed}
The set ${\cal N}(N)$  $\, (N>M)$ is imbedded in ${\cal N}(M)$,
i.e., ${\cal N}(N) \subset {\cal N}(M)$.
\end{Pn}
\begin{proof}.
By (\ref{1.10}) we have
\begin{equation}\label{2.58}
q(\lambda)^{N+1}W_{N+1}(\lambda)^*j W_{N+1}(\lambda)\leq
q(\lambda)^{M+1}W_{M+1}(\lambda)^*j W_{M+1}(\lambda), \quad
\lambda \in \BC_+.
\end{equation}
Insert the {\it length} $N$ of the interval into the notation
$\cal W$:
\begin{equation} \label{2.58'}
{\cal W}(N, \lambda)= {\cal W}( \lambda)=K W_{N+1}(
\ov{\lambda})^*.
\end{equation}
From (\ref{2.58}) and (\ref{2.58'}) it follows that
\begin{equation}\label{2.59}
q(\lambda)^{N-M}\Big({\cal W}(M, \lambda)^{-1}{\cal W}(N,
\lambda)\Big)^*j{\cal W}(M, \lambda)^{-1}{\cal W}(N, \lambda) \leq
j.
\end{equation}
Moreover, in view of (\ref{0.1}) and (\ref{2.58'}) we have
\begin{equation}\label{2.59'}
{\cal W}(M, \lambda)^{-1}{\cal W}(N,
\lambda)=\prod_{k=M+1}^{N}(I_m+\frac{i}{\lambda}C_k j),
\end{equation}
and the expression on the left-hand side of (\ref{2.59'}) is
analytic at $\lambda = -i$. Suppose now that $\vp \in {\cal N}(N)$
is a Weyl function generated by some pair $R$, $Q$, which
satisfies (\ref{2.55}). Then, according to (\ref{2.55}),
(\ref{2.59}) and (\ref{2.59'}) the pair
\begin{equation}\label{2.60}
\left[
\begin{array}{c}
\wt  R(\lambda) \\ \wt Q(\lambda)
\end{array}
\right]={\cal W}(M, \lambda)^{-1}{\cal W}(N, \lambda)\left[
\begin{array}{c}
 R(\lambda) \\  Q(\lambda)
\end{array}
\right]
\end{equation}
satisfies conditions of Definition \ref{genWeyl} too. Moreover, it
is easy to see that
\[
 i \bigl( {\cal W}_{21}(M, \lambda )\wt R(\lambda )+ {\cal W}_{22}(M, \lambda
)\wt Q(\lambda )\bigr) \bigl( {\cal W}_{11}(M, \lambda ) \wt
R(\lambda )+{\cal W}_{12}(M, \lambda )\wt Q(\lambda )\bigr)^{-1},
\]
\begin{equation}\label{2.61}
= i \bigl( {\cal W}_{21}(N, \lambda )R(\lambda )+ {\cal W}_{22}(N,
\lambda )Q(\lambda )\bigr) \bigl( {\cal W}_{11}(N, \lambda )
R(\lambda )+{\cal W}_{12}(N, \lambda )Q(\lambda )\bigr)^{-1}
\end{equation}
\[
=\vp(\lambda),
\]
which completes the proof.
\end{proof}

Theorem \ref{Tm2.2} and Proposition \ref{imbed} imply a
Borg-Marchenko type result:
\begin{Tm}\label{Tm2.4}
Let $\wt \varphi$ and $\wh \varphi$ be  Weyl functions of the two
discrete Dirac type systems (\ref{0.1}), which satisfy conditions
(\ref{2.4''}). Denote by $\wt C_k$ ($0 \leq k \leq \wt N$) the
potentials $C_k$ of  the first system and by $\wh C_k$ ($0 \leq k
\leq \wh N$) the potentials of the second system. Denote Taylor
coefficients of $i \wt \vp\Big(i \Big( \frac{z+1}{z-1} \Big)\Big)$
and $i \wh \vp\Big(i \Big( \frac{z+1}{z-1} \Big)\Big)$ at $z=0$ by
$\{\wt \a_k \}$ and $\{\wh \a_k \}$, respectively, and assume that
$\wt \a_k=\wh \a_k$ for all $k \leq N \leq \min\{ \wt N, \wh N\}$.
Then we have $\wt C_k=\wh C_k$ for $k \leq N$.
\end{Tm}
\begin{proof}.  According to Proposition \ref{imbed}, $\wt \varphi$ and $\wh \varphi$ are Weyl functions
of the first and second systems, respectively, on the interval $0
\leq k \leq N$. By Theorem \ref{Tm2.2} these systems on the
interval $0 \leq k \leq N$ are uniquely recovered by the first
$N+1$ Taylor coefficients of the Weyl functions.
\end{proof}
An interesting Borg-Marchenko type result for supersymmetric Dirac
difference operators have been obtained earlier in \cite{CGR}.

\section{Toeplitz matrices and Dirac system on the
semiaxis} \label{Toepl}
 \setcounter{equation}{0}
 By \cite{SaAul}, p. 116 it is easy to recover a block Toeplitz matrix $S$ which satisfies
(\ref{2.13}), where the blocks $\Phi_1$ and $\Phi_2$
 of $\Pi$ are   given by
(\ref{2.11}). Namely, we have
\begin{equation}\label{2.66}
S=\{s_{j-k} \}_{k,j=0}^N, \quad s_{-k}=\a_k=s_k^* \quad (k>0),
\quad s_0=s_0^*= \a_0+\a_0^*.
\end{equation}
Moreover, this $S$ is a unique solution of (\ref{2.13}). A
description of all extensions of $S$ preserving the number of
negative eigenvalues, which uses transfer matrix function $w_A$,
is given in \cite{SaAul} (see also Theorem 4.1 in \cite{SaAtepl})
in terms of the linear fractional transformation
\begin{equation}\label{2.67}
\wh \vp(\lambda)=\Big(\wh R(\lambda) w_{11}(\lambda)+\wh
Q(\lambda) w_{21}(\lambda)\Big)^{-1}\Big(\wh R(\lambda)
w_{12}(\lambda)+\wh Q(\lambda) w_{22}(\lambda)\Big),
\end{equation}
where $\{w_{kj}(\lambda) \}_{k,j=1}^2=w_A(N, \lambda)$, and the
meromorphic pairs $\wh R$, $\wh Q$ have $J$-property, i.e.,
\begin{equation}\label{2.68}
 \wh R(\lambda) \wh R(\lambda)^*+\wh Q(\lambda) \wh Q(\lambda)^*>0, \quad
\quad \wh R(\lambda) \wh Q(\lambda)^*+\wh Q(\lambda) \wh
R(\lambda)^* \geq 0, \quad \lambda \in \BC_+.
\end{equation}
In particular, for the case $S>0$, which is treated here, the
matrix functions $\displaystyle{\wh
\vp\Big(-\frac{i(z+1)}{2(z-1)}\Big)}$ are always analytic at $z=0$
and admit the Taylor representation
\begin{equation}\label{2.69}
\wh \vp\Big(-\frac{i(z+1)}{2(z-1)}\Big)=\wh s_0+\wh s_{-1}z+ \wh
s_{-2}z^2+ \ldots
\end{equation}
Our next statement follows from Theorem 4.1  \cite{SaAtepl}.
\begin{Tm}\label{ext}
Assume that $S=\{ s_{j-k} \}_{k,j=0}^N>0$, and fix $\a_0$ such
that $\a_0+\a_0^*=s_0$. Using (\ref{2.11}), (\ref{2.37}), and
(\ref{2.66}) introduce $\Pi =[\Phi_1 \quad \Phi_2]$ and
$\{w_{kj}(\lambda) \}_{k,j=1}^2=w_A(N, \lambda)$. Now, let matrix
functions $\wh \vp$ be given by (\ref{2.67}), where the pairs $\wh
R$, $\wh Q$ satisfy (\ref{2.68}) and are well defined at $\lambda
=\frac{i}{2}$. Then   the Taylor coefficients $\wh s_{-k}$  at
$z=0$ of the matrix functions $\displaystyle{\wh
\vp\Big(-\frac{i(z+1)}{2(z-1)}\Big)}$ satisfy relations
\begin{equation}\label{2.70}
\wh s_{-k}= s_{-k} \quad (0 < k \leq N), \quad \wh s_0= \a_0.
\end{equation}
Moreover, putting $s_{-k}= s_{-k}^*=\wh s_{-k}$ for $k>N$, we have
$\{ s_{j-k} \}_{k,j=0}^M\geq 0$ for all $M>N$. In other words, the
Taylor coefficients  of the matrix functions $\displaystyle{\wh
\vp\Big(-\frac{i(z+1)}{2(z-1)}\Big)}$ generate nonnegative
extensions of $S$. All the nonnegative extensions of $S$ are
generated in this way.
\end{Tm}
Taking into account that $w_A(N, \lambda)Jw_A(N,\ov \lambda)^*=J$,
we derive equality $\wh \vp = \wt \vp$ for the matrix function
\begin{equation}\label{2.67'}
\wt \vp(\lambda)=-\Big( w_{12}(\ov \lambda)^*\wt
R(\lambda)+w_{22}(\ov \lambda)^*\wt Q(\lambda)\Big)\Big(
w_{11}(\ov \lambda)^*\wt R(\lambda)+w_{21}(\ov \lambda)^*\wt
Q(\lambda)\Big)^{-1},
\end{equation}
where
\begin{equation}\label{2.68'}
 \wt R(\lambda)^* \wt R(\lambda)+\wt Q(\lambda)^* \wt Q(\lambda)>0, \quad
\quad \wh R(\lambda) \wt Q(\lambda)+\wh Q(\lambda) \wt R(\lambda)
= 0, \quad \lambda \in \BC_+.
\end{equation}
Notice also that relations (\ref{2.68}) and (\ref{2.68'}) yield
\begin{equation}\label{2.71}
\wt R(\lambda)^* \wt Q(\lambda)+\wt Q(\lambda)^* \wt R(\lambda)
\leq 0,
\end{equation}
and vice versa relations (\ref{2.68'}) and (\ref{2.71}) yield the
second relation in (\ref{2.68}). Hence, Theorem \ref{ext} can be
reformulated in terms of the linear fractional transformations
(\ref{2.67'}), where $\wt R$, $\wt Q$ have $J$-property
(\ref{2.71}). Finally, use (\ref{2.4}), (\ref{2.5'}) and
(\ref{2.36}) to rewrite (\ref{2.5}) in the form
\[
i \vp(\lambda)=-\Big( w_{12}(-\ov \lambda/2)^* \wt
R(-\lambda/2)+w_{22}(-\ov \lambda/2)^*\wt Q(-\lambda/2)\Big)
\]
\begin{equation}\label{2.72}
\times \Big( w_{11}(-\ov \lambda/2)^*\wt R(-\lambda/2)+w_{21}(-\ov
\lambda/2)^*\wt Q(-\lambda/2 )\Big)^{-1},
\end{equation}
where we put
\begin{equation}\label{2.73}
\left[
\begin{array}{c}
 \wt R(-\lambda/2) \\ \wt Q(-\lambda/2)
\end{array}
\right]=K\left[
\begin{array}{c}
 R(\lambda) \\  Q(\lambda)
\end{array}
\right].
\end{equation}
Here, formula (\ref{2.73}) is a one to one mapping of the pairs
satisfying (\ref{2.55}) into pairs satisfying the first relation
in (\ref{2.68'}) and relation (\ref{2.71}). By (\ref{2.67'}) and
(\ref{2.72}) we have $\displaystyle{\wh
\vp\Big(-\frac{i(z+1)}{2(z-1)}\Big)=\wt
\vp\Big(-\frac{i(z+1)}{2(z-1)}\Big)=i \vp
\Big(i\frac{(z+1)}{(z-1)}\Big)}$ Therefore Theorem \ref{ext} can
be rewritten.
\begin{Tm}\label{extn}
Assume that $S=\{ s_{j-k} \}_{k,j=0}^N>0$, fix $\a_0$ such that
$\a_0+\a_0^*=s_0$, and introduce $\cal W$ via (\ref{2.4}) and
(\ref{2.36}). Let matrix functions $ \vp$ be given by (\ref{2.5}),
where the pairs $R$, $Q$ satisfy (\ref{2.55}) and are well defined
at $\lambda =-i$. Then  $i\vp(-i)=\a_0$, and the following Taylor
coefficients $ \a_{k}$  at $z=0$ of the matrix functions
$\displaystyle{i \vp\Big(i\frac{(z+1)}{(z-1)}\Big)}$ satisfy
relations
\begin{equation}\label{2.74}
\a_{k}= s_{-k} \quad (0 < k \leq N).
\end{equation}
Moreover, putting $s_{-k}= s_{k}^*=\a_k$ for $k>N$, we have
$\{s_{j-k} \}_{k,j=0}^M\geq 0$ for all $M>N$. In other words, the
Taylor coefficients  of the matrix functions $\displaystyle{i
\vp\Big(i\frac{(z+1)}{(z-1)}\Big)}$ generate nonnegative
extensions of $S$. All the nonnegative extensions of $S$ are
generated in this way.
\end{Tm}
\begin{Rk}\label{imb}
By Definition \ref{genWeyl} and Theorem \ref{extn} the  Weyl
functions from the Weyl disk ${\cal N}(N)$ generate all the
nonnegative extentions of $S$. It provides, in particular,
another proof of Proposition \ref{imbed}.
\end{Rk}
\begin{equation}\label{2.75}
\end{equation}
Consider now system (\ref{0.1}), which satisfies (\ref{2.4''}) on
the semiaxis $k \geq 0$.
\begin{Tm}\label{Pn2.5} Let system (\ref{0.1})  be given on the semiaxis $k \geq 0$
and let matrices $C_k$ satisfy (\ref{2.4''}). Then, there is a
unique function $\vp_{\infty}$, which belongs to all the  Weyl
discs ${\cal N}(N)$:
\begin{equation}\label{2.62}
\bigcap_{N=0}^\infty {\cal N}(N)=\vp_{\infty}.
\end{equation}
\end{Tm}
\begin{proof}. According to Corollary \ref{CyTR} matrices $\{C_k\}_{k=0}^N$
($N < \infty$) or equivalently matrices $\{\b(k)\}_{k=0}^N$
uniquely define blocks $\{s_{-k}\}_{k=0}^N$, where
$s_{-k}=\a_k=\psi_k-\psi_{k-1}$ for $k>0$ and $s_0=\a_0+\a_0^* \,$
($\a_0 = \psi_0$). Moreover, by Proposition \ref{imbed} these
$s_{-k}$ do not depend on $N \geq k$, and so system (\ref{0.1})
on the semiaxis determines an infinite sequence
$\{s_{-k}\}_{k=0}^\infty$. By  Theorem \ref{Tm2.2} we have $\{
s_{j-k} \}_{k,j=0}^N>0$ for all $N \geq 0$. Apply now Theorem
\ref{extn} to see that
\begin{equation}\label{2.63}
i\vp\Big(i\frac{(z+1)}{(z-1)}\Big)=\a_0+\sum_{k=1}^{\infty}s_{-k}z^k=
i\vp_{\infty}\Big(i\frac{(z+1)}{(z-1)}\Big),
\end{equation}
i.e., this $\vp$ belongs to $\bigcap_{N=0}^\infty {\cal N}(N)$.
Moreover, as the sequence $\{s_{-k}\}_{k=0}^\infty$ is unique, so
by Theorem \ref{extn} the function $\vp \in \bigcap_{N=0}^\infty
{\cal N}(N)$ is unique.
\end{proof}

Recall that a Weyl function on the semiaxis is defined by
Definition \ref{WT}, where $K$ is given by formula (\ref{1.26}).
Theorem \ref{Pn2.5} yields our next result.
\begin{Tm} \label{Wonsa}
Let system (\ref{0.1})  be given on the semiaxis $k \geq 0$ and
let matrices $C_k$ satisfy (\ref{2.4''}). Then, the matrix
function $\vp_{\infty}$ given by (\ref{2.62})  is the unique Weyl
function of system (\ref{0.1}) on the semiaxis.
\end{Tm}
\begin{proof}.
By (\ref{2.5}) and (\ref{2.62}),  we have
\begin{equation} \label{2.77}
\left[\begin{array}{c}
 -i\vp_{\infty}(\lambda) \\ I_p
\end{array}
\right]=J {\cal W}(r+1, \lambda) \left[\begin{array}{c}
 R(\lambda) \\ Q(\lambda)
\end{array}
\right]
\end{equation}
for all $r\geq 0$ and for some depending on $r$ pairs $R$, $Q$,
which satisfy (\ref{2.55}). In view of (\ref{2.4}), (\ref{2.5'}),
(\ref{2.45}), and (\ref{2.77}) we obtain
\[
[i\vp_{\infty}(\lambda)^* \quad I_p ]\Big(q(\lambda)^{r+1} K
W_{r+1}(\lambda)^*jW_{r+1}(\lambda)K^*-J\Big)\left[\begin{array}{c}
 -i\vp_{\infty}(\lambda) \\ I_p
\end{array}
\right]
\]
\begin{equation} \label{2.78}
=i\Big(\vp_{\infty}(\lambda)-\vp_{\infty}(\lambda)^*\Big)+\left(\frac{|\lambda^2+1|^2}{|\lambda|^2(|\lambda|^2+1)}
\right)^{r+1}[R(\lambda)^* \quad
Q(\lambda)^*]j\left[\begin{array}{c}
 R(\lambda) \\ Q(\lambda)
\end{array}
\right].
\end{equation}
Now, formulas (\ref{2.55}) and (\ref{2.78}) imply
\[
[i\vp_{\infty}(\lambda)^* \quad I_p ]\Big(q(\lambda)^{r+1} K
W_{r+1}(\lambda)^*jW_{r+1}(\lambda)K^*-J\Big)\left[\begin{array}{c}
 -i\vp_{\infty}(\lambda) \\ I_p
\end{array}
\right]
\]
\begin{equation} \label{2.79}
\leq i\Big(\vp_{\infty}(\lambda)-\vp_{\infty}(\lambda)^*\Big).
\end{equation}
It follows from (\ref{1.10}) and (\ref{2.79}) that
\[
\sum_{k=0}^r[i\vp_{\infty}(\lambda)^* \quad I_p ]q(\lambda)^k K
W_k(\lambda)^*C_kW_k(\lambda)K^*\left[\begin{array}{c}
 -i\vp_{\infty}(\lambda) \\ I_p
\end{array}
\right]
\]
\begin{equation} \label{2.80}
 \leq \frac{|\lambda|^2+1}{\lambda - \ov \lambda} \Big(\vp_{\infty}(\lambda)-\vp_{\infty}(\lambda)^*\Big).
\end{equation}
Finally, by (\ref{2.80}) the inequality (\ref{2.76}) is immediate,
and $\vp_{\infty}$ given by (\ref{2.62})  is a Weyl function.

To prove the uniqueness of the Weyl function notice that by
Proposition \ref{Pn3.2} and by relation (\ref{1.11}) we have
\begin{equation} \label{2.81}
q^k W_k^*C_kW_k \geq q^k  W_k^*jW_k \geq q^{k-1}W_{k-1}^*jW_{k-1}
\geq \ldots \geq W_{0}^*jW_{0}=j \quad (\lambda \in \BC_-).
\end{equation}
Hence, in view of (\ref{2.81}) we obtain
\[
\sum_{k=0}^r [I_p \quad I_p ]q(\lambda)^k K
W_k(\lambda)^*C_kW_k(\lambda)K^*\left[\begin{array}{c}
 I_p \\ I_p
\end{array}
\right] \geq  2(r+1)I_p,
\]
and it follows that
\begin{equation} \label{2.82}
\sum_{k=0}^\infty [I_p \quad I_p ]q(\lambda)^k K
W_k(\lambda)^*C_kW_k(\lambda)K^*\left[\begin{array}{c}
 I_p \\ I_p
\end{array}
\right] = \infty .
\end{equation}
Taking into account Definition \ref{WT} and inequality
(\ref{2.82}),  we can show the uniqueness of the Weyl function
similar to the proof of the  uniqueness in Theorem \ref{TmDirect}.
\end{proof}
Now, we formulate a solution of the inverse problem.
\begin{Tm} \label{genInv}
The set of the Weyl functions $\vp(\lambda)$ of systems
(\ref{0.1}), given on the semiaxis $k \geq 0$ and such that
matrices $C_k$ satisfy (\ref{2.4''}), coincides with the set of
functions $\vp$ such that
\begin{equation}\label{2.85}
i\vp\Big(i\frac{(z+1)}{(z-1)}\Big)=\a_0+\sum_{k=1}^{\infty}s_{-k}z^k
\end{equation}
are Caratheodory matrix functions in the unit disk and $\{ s_{j-k}
\}_{k,j=0}^N>0$ for all $0 \leq N<\infty$ $\,(s_0:=\a_0+\a_0^*)$.
These systems (\ref{0.1}) are uniquely recovered from their Weyl
functions via the procedure given in Theorem \ref{Tm2.2}.
\end{Tm}
\begin{proof}. According to Theorem  \ref{Wonsa} the Weyl function on the semiaxis
is also a Weyl function on the intervals. Hence, the procedure to
construct a solution of the inverse problem follows from Theorem
\ref{Tm2.2}. It follows from Theorem \ref{Tm2.2} also, that the
matrices $\{ s_{j-k} \}_{k,j=0}^N$ generated by the Weyl functions
are positive definite.

Hence, it remaines to show that all the functions such that
(\ref{2.85}) holds and $\{ s_{j-k} \}_{k,j=0}^r>0$ ($r \geq 0$)
are Weyl functions. Indeed, fixing such a matrix function $\vp$,
we get a sequence of matrices $S(r)=\{ s_{j-k} \}_{k,j=0}^r>0$.
Therefore we get a sequence of the transfer matrix functions
$w_A(r, \lambda)$ of the form (\ref{2.37}), where
\begin{equation}\label{2.88}
\Pi(r)=\left[
\begin{array}{ll}
I_{p} & \a_0 \\ I_{p} & \a_0+s_{-1} \\ \ldots & \ldots\\I_p &
\a_0+s_{-1}+ \ldots+s_{-r}
\end{array}
\right],
\end{equation}
and  (\ref{2.19}) holds. Taking into account formulas
(\ref{2.37}), (\ref{2.43}) and (\ref{2.39}), (\ref{2.42}), we
obtain matrix functions
\begin{equation}\label{2.86}
\b(r)^*\b(r):= \Pi(r)^*S(r)^{-1}P^*\big(PS(r)^{-1}P^*\big)^{-1} P
S(r)^{-1}\Pi(r) \quad (r \geq 0).
\end{equation}
In view of the matrix identity (\ref{2.19})   we have
\[
 P S(r)^{-1}\Pi(r)J \Pi(r)^*S(r)^{-1}P^*=-iP
\Big(S(r)^{-1}A(r)-A(r)^*S(r)^{-1}\Big)P^*
\]
\begin{equation}\label{2.87}
=PS(r)^{-1}P^*.
\end{equation}
In other words $\b(r)$ satisfies the second relation in
(\ref{2.4''}). Therefore formulas $C_r=2K \b(r)^*\b(r)K-j$ define
a system of our class on the semiaxis. In fact, this construction
coincides with  (\ref{2.14}) and the function $\vp$ is the Weyl
function of our system. Indeed, similar to the proof of Theorem
\ref{Tm2.2} we derive from (\ref{2.86})  the equality
(\ref{2.36}). Compare now  Definition \ref{genWeyl} and Theorem
\ref{extn} to see that $\vp \in {\cal N}(N)$ for any $N$.
According to Theorems \ref{Pn2.5} and \ref{Wonsa} it means that
$\vp$ is the Weyl function.
\end{proof}

Finally, consider  the upper halfplane and define holomorphic Weyl
functions in $\BC_+$ via relations (\ref{2.76}) and (\ref{1.26})
too. Put
\begin{equation} \label{2.83}
\ov{\cal N}(N):=\{\vp(\ov \lambda)^*: \, \vp \in {\cal N}(N) \}.
\end{equation}
\begin{Rk} \label{upper}
Similar to the proof that $\wh \vp=\wt \vp$, where  $\wh \vp$ and
$\wt \vp$ are given by (\ref{2.67}) and (\ref{2.67'}),
respectively, one can show that the set $\ov{\cal N}(N)$ consists
of linear fractional transformations (\ref{2.5}), where the pairs
$R$, $Q$ are meromorphic in $\BC_+$, are well defined at $\lambda
= i$, and have the property
\begin{equation}\label{2.84}
 R(\lambda)^*  R(\lambda)+ Q(\lambda)^* Q(\lambda)>0, \quad
\quad R(\lambda)^* R(\lambda)\geq  Q(\lambda)^* Q(\lambda), \quad
\lambda \in \BC_+.
\end{equation}
\end{Rk}
In view of Remark \ref{upper}, we obtain in $\BC_+$ the analog of
Theorem \ref{Wonsa}, and the proof is similar.
\begin{Tm} \label{Wons}
Let system (\ref{0.1})  be given on the semiaxis $k \geq 0$ and
let matrices $C_k$ satisfy (\ref{2.4''}). Then, the matrix
function $\vp_{\infty}(\ov \lambda)^*=\bigcap_{N=0}^\infty \ov
{\cal N}(N)$   is the unique Weyl function in $\BC_+$ of system
(\ref{0.1}) on the semiaxis.
\end{Tm}
\begin{proof}.
Substitute $\vp_{\infty}(\ov \lambda)^*$ instead of $\vp_{\infty}(
\lambda)$ into (\ref{2.78}) and take into account (\ref{2.84}) to
derive
\[
[i\vp_{\infty}(\ov \lambda) \quad I_p ]\Big(J-q(\lambda)^{r+1} K
W_{r+1}(\lambda)^*jW_{r+1}(\lambda)K^*\Big)\left[\begin{array}{c}
 -i\vp_{\infty}(\ov \lambda)^* \\ I_p
\end{array}
\right]
\]
\begin{equation} \label{2.79'}
\leq i\Big(\vp_{\infty}(\ov \lambda)-\vp_{\infty}(\ov
\lambda)^*\Big).
\end{equation}
Now, inequalities (\ref{2.80}) and  (\ref{2.76}) are
straightforward, i.e.,
 $\vp_{\infty}(\ov \lambda)^*$ is a Weyl function.

Instead of (\ref{2.81}) we use the inequality
\begin{equation} \label{2.81'}
q(\lambda)^k W_k(\lambda)^*C_kW_k(\lambda) \geq -q(\lambda)^k
W_k(\lambda)^*jW_k(\lambda) \geq -j \quad (\lambda \in \BC_+),
\end{equation}
which yields inequality
\begin{equation} \label{2.82'}
\sum_{k=0}^\infty [I_p \quad -I_p ]q(\lambda)^k K
W_k(\lambda)^*C_kW_k(\lambda)K^*\left[\begin{array}{c}
 I_p \\ -I_p
\end{array}
\right] = \infty .
\end{equation}
The uniqueness of the Weyl function follows from (\ref{2.82'})
\end{proof}
Theorem \ref{Wons} for the  scalar case $p=1$ has been proved
earlier in \cite{Ger, GoNe} (see also Theorem 3.2.11 \cite{Si2}).

{\bf Acknowledgements.} The research of I. Roitberg was supported
by Otto Benecke Stiftung. The work of A.L. Sakhnovich was
supported by the Austrian Science Fund (FWF) under Grant  no.
Y330.

B. Fritzsche,  Fakult\"at f\"ur Mathematik und Informatik, \\
Mathematisches Institut, Universit\"at Leipzig, Augustusplatz
10/11, \\ D-04109 Leipzig, Germany, fritzsche@math.uni-leipzig.de

B. Kirstein, Fakult\"at f\"ur Mathematik und Informatik, \\
Mathematisches Institut, Universit\"at Leipzig, Augustusplatz
10/11, \\ D-04109 Leipzig, Germany, kirstein@math.uni-leipzig.de

I.Ya. Roitberg, Fakult\"at f\"ur Mathematik und Informatik, \\
Mathematisches Institut, Universit\"at Leipzig, Augustusplatz
10/11, \\ D-04109 Leipzig, Germany, i$_-$roitberg@yahoo.com

A.L. Sakhnovich, Fakult\"at f\"ur Mathematik, Universit\"at Wien,
\\
 Nordbergstrasse 15, A-1090 Wien, Austria,
 al$_-$sakhnov@yahoo.com

\end{document}